\documentclass[12pt]{article}%
\usepackage{amsmath}
\usepackage{amsfonts}
\usepackage{amssymb}
\usepackage{graphicx}
\usepackage{color}%
\providecommand{\U}[1]{\protect\rule{.1in}{.1in}}
\newtheorem{theorem}{Theorem}[section]

\newtheorem{corollary}[theorem]{Corollary}

\newtheorem{definition}[theorem]{Definition}

\newtheorem{lemma}[theorem]{Lemma}
\newtheorem{notation}[theorem]{Notation}

\newtheorem{proposition}[theorem]{Proposition}
\newtheorem{remark}[theorem]{Remark}

\newcommand{\BIGOP}[1]{\mathop{\mathchoice{\raise-0.22em\hbox{\huge
$#1$}} {\raise-0.05em\hbox{\Large $#1$}}{\hbox{\large $#1$}}{#1}}}
\newenvironment{proof}[1][Proof]{\textbf{#1.} }{\ \rule{0.5em}{0.5em}}

\makeatletter\@addtoreset{equation}{section}\makeatother
\evensidemargin=\oddsidemargin
\newdimen\dummy
\dummy=\oddsidemargin
\begin{document}

\title{The Green's function for an acoustic, half-space impedance problem\\Part II: Analysis of the slowly varying and the plane wave component.}
\author{C. Lin\thanks{(chuhe.lin@math.uzh.ch), Institut f\"{u}r Mathematik,
Universit\"{a}t Z\"{u}rich, Winterthurerstr 190, CH-8057 Z\"{u}rich,
Switzerland}
\and J.M. Melenk\thanks{(melenk@tuwien.ac.at), Institut f\"{u}r Analysis und
Scientific Computing, Technische Universit\"{a}t Wien, Wiedner Hauptstrasse
8-10, A-1040 Wien, Austria.}
\and S. Sauter\thanks{(stas@math.uzh.ch), Institut f\"{u}r Mathematik,
Universit\"{a}t Z\"{u}rich, Winterthurerstr 190, CH-8057 Z\"{u}rich,
Switzerland}}
\maketitle

\begin{abstract}
We show that the acoustic Green's function for a half-space impedance problem
in arbitrary spatial dimension $d$ can be written as a sum of two terms, each
of which is the product of an exponential function with the \textit{eikonal}
in the argument and a \textit{slowly varying} function. We introduce the
notion of families of slowly varying functions to formulate this statement as
a theorem and present its proof.

\end{abstract}



\section{Introduction}

Wave phenomena in half-space domains have many important applications such as
the modelling of noise propagation over flat ground, the detection of
inclusions in homogeneous soil under a plane crust, as well as the design of
noise barriers positioned over a large flat boundary of an ambient domain
modelled by a half-space. If the problem is formulated in the frequency domain
and modelled by a Helmholtz equation, the boundary element method is a
well-established and popular numerical discretization method for acoustic
problems in unbounded domains. Typically it involves the explicit knowledge of
the full-space Green's function of the underlying differential operators. For
half-space problems, however, the Green's function for the half-space is
involved -- for sound-hard and sound-soft boundary conditions this function is
well-known and can be explicitly expressed by the classical method of images.
For the physically more relevant case of impedance boundary conditions, the
Green's function in general is not known explicitly and typically expressed by
oscillatory Fourier-type integrals, Hankel transforms, and Sommerfeld-type
integrals (see, e.g., \cite[(13)]{ChandlerWilde_imp_Green}, \cite[(21)]%
{Hein_Nedelec_green}, \cite{Duran_stat_phase}, \cite{Duran_stat_phase_long},
\cite{hoernig2010green}, \cite{ochmann2011closed}, \cite{Gimperlein_halfspace}%
, \cite{ROJAS20183903}). While standard methods from asymptotic analysis such
as the method of stationary phase allow for detecting the lowest order term in
an asymptotic expansion it is a fairly open problem to obtain uniform higher
order asymptotic expansions for this Green's function as well as estimates of
the remainder and structural insight in this function.

It is the main goal of this paper to report on progress in the described
direction. The starting point is a new representation of the half-space
acoustic impedance Green's function in arbitrary spatial dimension
$d\in\left\{  1,2,3,\ldots\right\}  $ that is derived in the companion paper
\cite{LinMelenkSauter_Gimp_I}. In contrast to the representations via
oscillatory integrals, the integrand in the new integral representation is
non-oscillatory with respect to the outer variable and defines a function that
is non-oscillatory in a sense which will be made precise in this paper.

The analysis in our paper is based on this new representation and makes the
following statement rigorous:%
\begin{equation}%
\begin{array}
[c]{l}%
\text{The acoustic Green's function for a half-space impedance problem can be
written as the}\\
\text{sum of (two) terms, each of which is the product of an oscillatory
exponential function}\\
\text{and a slowly varying one.}%
\end{array}
\label{softstatement}%
\end{equation}
This is in the spirit of the theory of geometric optics and ray theory;
standard references include \cite[Chap.~2]{Babich_book}, \cite[Chap.~3]%
{borovikov1994geometrical}, \cite{Buslaev_1975} where scattering problems are
analysed as (the sum of) products of the exponential function with the eikonal
in the argument and a slowly varying function. Here our focus is the
development of a rigorous analysis of the slowly varying function for the
specific half-space problem.

This analysis has immediate implication for its approximation: the slowly
varying part can be approximated by polynomials with exponential convergence
and the oscillatory exponential factor can be treated by directional
hierarchical ($\mathcal{DH}^{2}$) matrices (see \cite{Engquist_Ying,
Boerm_Melenk,brm2019variable}) or by a butterfly algorithm
\cite{Demanet_butterfly1}.

The paper is structured as follows. In Section~\ref{SecGreenFct} we formulate
the governing equations for the acoustic half-space Green's function and
briefly recall the representation from \cite{LinMelenkSauter_Gimp_I}. In
Section~\ref{SecAnaGreen} we develop the analysis of the Green's function so
that statement (\ref{softstatement}) can be made rigorous in Theorem
\ref{TheoSlowVar}. This requires as a prerequisite the notion of slowly
varying families of functions (Sec.~\ref{SecSlowVaryFct}), the derivation of a
majorant of the Bessel function of second kind which is immanent in the
representation of the Green's function (Sec.~\ref{SecMacdonald1}), the
investigation of holomorphic norm extensions (Sec.~\ref{SecHoloNorm}), and
finally the analysis of a coordinate transform which is involved in the
integral representation of the impedance part of the half-space Green's
function (see Sec.~\ref{Secmue}).

\section{The Green's function for the acoustic half-space problem with
impedance boundary conditions\label{SecGreenFct}}

Let the upper/lower half-space in $\mathbb{R}^{d}$, $d\in\left\{
1,2,\ldots\right\}  $, and its boundary be denoted by%
\begin{align*}
H_{+}  &  :=\left\{  \mathbf{x}=\left(  x_{j}\right)  _{j=1}^{d}\in
\mathbb{R}^{d}\mid x_{d}>0\right\}  ,\\
H_{-}  &  :=\left\{  \mathbf{x}=\left(  x_{j}\right)  _{j=1}^{d}\in
\mathbb{R}^{d}\mid x_{d}<0\right\}  ,\\
H_{0}  &  :=\partial H_{+}:=\left\{  \mathbf{x}=\left(  x_{j}\right)
_{j=1}^{d}\in\mathbb{R}^{d}\mid x_{d}=0\right\}
\end{align*}
with outward normal vector $\mathbf{n}=\left(  0,\ldots,0,-1\right)  ^{T}$.
Let
\[
\overset{\bullet}{\mathbb{C}}_{\geq0}:=\left\{  \zeta\in\mathbb{C}%
\mid\operatorname{Re}\zeta\geq0\right\}  \backslash\left\{  0\right\}
\quad\text{and\quad}\mathbb{C}_{>0}:=\left\{  \zeta\in\mathbb{C}%
\mid\operatorname{Re}\zeta>0\right\}  .
\]

We consider the problem to find the Green's function $G:H_{+}\times
H_{+}\rightarrow\mathbb{C}$ for the acoustic half-plane problem with impedance
boundary conditions:%
\begin{equation}%
\begin{array}
[c]{rll}%
-\Delta_{\mathbf{x}}G\left(  \mathbf{x},\mathbf{y}\right)  +s^{2}G\left(
\mathbf{x},\mathbf{y}\right)  & =\delta_{0}\left(  \mathbf{x}-\mathbf{y}%
\right)  & \text{for }\left(  \mathbf{x},\mathbf{y}\right)  \in H_{+}\times
H_{+},\\
\frac{\partial}{\partial\mathbf{n}_{\mathbf{x}}}G\left(  \mathbf{x}%
,\mathbf{y}\right)  +s\beta G\left(  \mathbf{x},\mathbf{y}\right)  & =0 &
\text{for }\left(  \mathbf{x},\mathbf{y}\right)  \in H_{0}\times H_{+}\\
G\left(  r%
\mbox{\boldmath$ \zeta$}%
,\mathbf{y}\right)  & \overset{r\rightarrow+\infty}{\rightarrow}0 & \text{for
}\left(
\mbox{\boldmath$ \zeta$}%
,\mathbf{y}\right)  \in H_{+}\times H_{+}.
\end{array}
\label{goveq}%
\end{equation}
for some $\beta>0$ and frequency $s\in\mathbb{C}_{>0}$. The index $\mathbf{x}$
in the differential operators indicates that differentiation is done with
respect to the variable $\mathbf{x}$.

\begin{remark}
\label{RemLimitAbsorption}Problem (\ref{goveq}) is formulated for
$s\in\mathbb{C}_{>0}$. The Green's function $G=G_{s}$ depends on $s$ and for
$\operatorname{Re}s>0$ it is assumed to decay for $\mathbf{x}=r%
\mbox{\boldmath$ \zeta$}%
$ as $r\rightarrow+\infty$ for any fixed direction $%
\mbox{\boldmath$ \zeta$}%
\in H_{+}$. Problem (\ref{goveq}) for the case $s\in\operatorname*{i}%
\mathbb{R}\backslash\left\{  0\right\}  $ is considered as the limit from the
positive complex half-plane $\mathbb{C}_{>0}$:%
\[
G_{s}=\lim_{\substack{\zeta\rightarrow s\\\zeta\in\mathbb{C}_{>0}}}G_{\zeta
}\text{.}%
\]

\end{remark}

In the following, we focus our attention on spatial dimensions $d\geq2$ and
general impedance parameter $\beta>0$. In \cite{LinMelenkSauter_Gimp_I}, fully
explicit representations of the Green's function are presented for $d=1$ and
for $d>1$ with $\beta=1$. Next we recall the new integral representation given
in \cite{LinMelenkSauter_Gimp_I}.

The representation of the Green's function as the solution of (\ref{goveq})
requires some preparations. Let $K_{\nu}$ denote the Macdonald function
(modified Bessel function of second kind and order $\nu$, see, e.g.,
\cite[\S 10.25]{NIST:DLMF}, \cite{Macdonald_fct_orginal}). We introduce the
function%
\begin{equation}
g_{\nu}\left(  r\right)  :=\frac{1}{\left(  2\pi\right)  ^{\nu+3/2}}\left(
\frac{s}{r}\right)  ^{\nu+1/2}K_{\nu+1/2}\left(  sr\right)  \label{defgnue}%
\end{equation}
and note that $g_{\nu}\left(  \left\Vert \mathbf{x}-\mathbf{y}\right\Vert
\right)  $ is the full space Green's function for the Helmholtz operator (see
\cite[Thm. 4.4]{Mclean00} and \cite[(6), (12)]{Buslaev_1975} in combination
with the connecting formula \cite[\S 10.27.8]{NIST:DLMF}). For $\mathbf{y=}%
\left(  y_{j}\right)  _{j=1}^{d}\in H_{+}$, we introduce the reflection
operator $\mathbf{Ry}=\left(  \mathbf{y}^{\prime},-y_{d}\right)  $, where
$\mathbf{y}^{\prime}=\left(  y_{j}\right)  _{j=1}^{d-1}$. The dependence on
the spatial dimension $d$ will be expressed via the parameter%
\[
\nu:=\left(  d-3\right)  /2.
\]
Let%
\[
\mathbb{Z}^{\operatorname*{half}}:=\left\{  -\frac{1}{2},0,\frac{1}%
{2},1,\ldots\right\}  \quad\text{and for }\mu\geq-\frac{1}{2}:\text{\quad
}\mathbb{Z}_{\geq\mu}^{\operatorname*{half}}:=\left\{  \nu\in\mathbb{Z}%
^{\operatorname*{half}}\mid\nu\geq\mu\right\}  .
\]
so that $d\in\left\{  2,3,\ldots\right\}  $ is equivalent to $\nu\in
\mathbb{Z}^{\operatorname*{half}}$.

Let the functions $r:\mathbb{R}^{d}\rightarrow\mathbb{R}$ and $r_{+}%
:\mathbb{R}^{d}\rightarrow\mathbb{R}$ be defined for $\mathbf{z}\in H_{+}$ and
$\mathbf{z}^{\prime}:=\left(  z_{j}\right)  _{j=1}^{d-1}$ by%

\begin{equation}
r\left(  \mathbf{z}\right)  :=\left\Vert \mathbf{z}\right\Vert ,\quad
r_{+}\left(  \mathbf{z}\right)  :=r\left(  \mathbf{z}\right)  +\beta z_{d}
\label{Defrrplus}%
\end{equation}
and set%
\begin{equation}
y\left(  \mathbf{z,\cdot}\right)  :\left[  z_{d},\infty\right[  \rightarrow
\left[  0,\infty\right[  , \qquad y\left(  \mathbf{z},t\right)  :=-r_{+}%
\left(  \mathbf{z}\right)  +\beta t+\mu\left(  \mathbf{z}^{\prime},t\right)
\label{defyvont}%
\end{equation}
with the function $\mu\left(  \mathbf{z}^{\prime},\cdot\right)  :\left[
z_{d},\infty\right[  \rightarrow\left[  \left\Vert \mathbf{z}\right\Vert
,\infty\right[  $ given by%
\[
\mu\left(  \mathbf{z}^{\prime},t\right)  :=\sqrt{\left\Vert \mathbf{z}%
^{\prime}\right\Vert ^{2}+t^{2}}.
\]
The derivative of $y$ satisfies%
\begin{equation}
\frac{\partial y\left(  \mathbf{z},t\right)  }{\partial t}=\beta+\frac{t}%
{\mu\left(  \mathbf{z}^{\prime},t\right)  }>0 \label{cfdydt}%
\end{equation}
so that $y\left(  \mathbf{z},\cdot\right)  $ maps the interval $\left[
z_{d},\infty\right[  $ strictly increasing onto $\left[  0,\infty\right[  $.
Its inverse%
\begin{equation}
t\left(  \mathbf{z},\cdot\right)  :\left[  0,\infty\right[  \rightarrow\left[
z_{d},\infty\right[  \label{tsub}%
\end{equation}
is also strictly increasing. The derivative $\partial t\left(  \mathbf{z}%
,y\right)  /\partial y$ can be expressed by using (\ref{cfdydt}):%
\begin{equation}
\frac{\partial t\left(  \mathbf{z},y\right)  }{\partial y}=\frac{\tilde{\mu
}\left(  \mathbf{z},y\right)  }{t\left(  \mathbf{z},y\right)  +\beta\tilde
{\mu}\left(  \mathbf{z},y\right)  }, \label{dert}%
\end{equation}
where%
\begin{equation}
\tilde{\mu}\left(  \mathbf{z},y\right)  :=\mu\left(  \mathbf{z}^{\prime
},t\left(  \mathbf{z},y\right)  \right)  \text{\quad and\quad}\frac
{\partial\tilde{\mu}\left(  \mathbf{z},y\right)  }{\partial y}=\frac{t\left(
\mathbf{z},y\right)  }{t\left(  \mathbf{z},y\right)  +\beta\tilde{\mu}\left(
\mathbf{z},y\right)  }>0. \label{defsmue}%
\end{equation}
In the following, the shorthands%
\begin{equation}
t=t\left(  \mathbf{z},y\right)  ,\quad\tilde{\mu}=\tilde{\mu}\left(
\mathbf{z},y\right)  \label{shorthands}%
\end{equation}
will be used. A key role for the representation of the Green's function will
be played by the functions%
\begin{equation}
\sigma_{\nu}\left(  r,z\right)  :=\frac{z-\beta r}{z+\beta r}g_{\nu}\left(
r\right)  \quad\text{and\quad}\psi_{\nu,s}\left(  \mathbf{z}\right)
:=\frac{1}{s}\int_{0}^{\infty}\operatorname*{e}\nolimits^{-sy}q_{\nu}\left(
\mathbf{z},y\right)  dy \label{Defpsinues}%
\end{equation}
with%
\begin{equation}
q_{\nu}\left(  \mathbf{z},y\right)  :=\frac{d}{dy}\left(  \frac
{\operatorname*{e}\nolimits^{s\tilde{\mu}}K_{\nu+1/2}\left(  s\tilde{\mu
}\right)  }{\left(  t+\beta\tilde{\mu}\right)  \left(  s\tilde{\mu}\right)
^{\nu-1/2}}\right)  . \label{defqnue}%
\end{equation}
From \cite[Thm.~{3.1} and Rem.~{3.2}]{LinMelenkSauter_Gimp_I} the following
representation of the half-space Green's function follows.

\begin{definition}
Let $d\in\left\{  2,\ldots\right\}  $ denote the spatial dimension. The
Green's function for the acoustic half-space problem with impedance boundary
conditions, i.e., the solution of (\ref{goveq}) is given by%
\begin{equation}
G_{\operatorname*{half}}\left(  \mathbf{x},\mathbf{y}\right)
:=G_{\operatorname*{illu}}\left(  \mathbf{x}-\mathbf{y}\right)
+G_{\operatorname*{refl}}\left(  \mathbf{x}-\mathbf{Ry}\right)
+G_{\operatorname*{imp}}\left(  \mathbf{x}-\mathbf{Ry}\right)  ,
\label{Ghalfdef}%
\end{equation}
where $\nu=\left(  d-3\right)  /2$ and%
\[
G_{\operatorname*{illu}}\left(  \mathbf{z}\right)  :=g_{\nu}\left(  \left\Vert
\mathbf{z}\right\Vert \right)  ,\quad G_{\operatorname*{refl}}\left(
\mathbf{z}\right)  :=\sigma_{\nu}\left(  \left\Vert \mathbf{z}\right\Vert
,z_{d}\right)  ,\quad G_{\operatorname*{imp}}\left(  \mathbf{z}\right)
:=-\frac{\beta}{\pi}\left(  \frac{s^{2}}{2\pi}\right)  ^{\nu+1/2}%
\operatorname*{e}\nolimits^{-s\left\Vert \mathbf{z}\right\Vert }\psi_{\nu
,s}\left(  \mathbf{z}\right)
\]

\end{definition}

In \cite{LinMelenkSauter_Gimp_I} it proved that $G_{\operatorname*{half}}$
satisfies problem (\ref{goveq}), and thus the name Green's function is justified.

\section{Analysis of the half-space Green's function\label{SecAnaGreen}}

In this section, we will prove that the half-space Green's function in
(\ref{Ghalfdef}) can be split into a sum of (two) terms each of which can be
written as a product of the form $\operatorname*{e}^{-s\tau_{\ell}\left(
\mathbf{x},\mathbf{y}\right)  }g_{\ell}\left(  \mathbf{x},\mathbf{y}\right)
$, where $\tau_{\ell}$ denotes an \textit{eikonal} and the function $g_{\ell}$
is \textit{slowly varying.} To make the meaning of these notions precise some
preparations are necessary. For two points $\mathbf{x},\mathbf{y}$ in the
domain, the eikonal $\tau_{\ell}\left(  \mathbf{x},\mathbf{y}\right)  $ is the
length of a possible path of light (in the limit of geometric optics) for a
ray emitted at $\mathbf{y}$ and received at $\mathbf{x}$. In this way, the
eikonal for the direct ray between $\mathbf{x}$ and $\mathbf{y}$ is given by
$\tau_{\operatorname*{illu}}\left(  \mathbf{x},\mathbf{y}\right)  :=\left\Vert
\mathbf{x}-\mathbf{y}\right\Vert $ and for the reflected ray by $\tau
_{\operatorname*{refl}}\left(  \mathbf{x},\mathbf{y}\right)  :=\left\Vert
\mathbf{x}-\mathbf{Ry}\right\Vert $. For the half-space Green's function, we
set for $\mathbf{x},\mathbf{y}\in H_{+}$,
\begin{equation}
\mathbf{z}:=\mathbf{x}-\mathbf{y},\quad r:=r\left(  \mathbf{x},\mathbf{y}%
\right)  :=\left\Vert \mathbf{z}\right\Vert ,\quad\mathbf{z}_{-}%
:=\mathbf{x}-\mathbf{Ry}\text{,\quad}r_{-}:=r\left(  \mathbf{x},\mathbf{Ry}%
\right)  :=\left\Vert \mathbf{z}_{-}\right\Vert \label{defn}%
\end{equation}
and define%
\begin{subequations}
\label{slowGreen}
\end{subequations}%
\begin{align}
\Theta_{\nu,s}^{\operatorname*{illu}}\left(  \mathbf{x},\mathbf{y}\right)   &
:=\operatorname*{e}\nolimits^{sr}g_{\nu}\left(  r\right)  ,\tag{%
\ref{slowGreen}%
a}\label{slowGreena}\\
\Theta_{\nu,s}^{\operatorname*{refl}}\left(  \mathbf{x},\mathbf{y}\right)   &
:=\operatorname*{e}\nolimits^{sr_{-}}\sigma_{\nu}\left(  r_{-},x_{d}%
+y_{d}\right)  ,\tag{%
\ref{slowGreen}%
b}\label{slowGreenb}\\
\Theta_{\nu,s}^{\operatorname*{imp}}\left(  \mathbf{x},\mathbf{y}\right)   &
:=-\frac{\beta}{\pi}\left(  \frac{s^{2}}{2\pi}\right)  ^{\nu+1/2}\psi_{\nu
,s}\left(  \mathbf{z}_{-}\right)  . \tag{%
\ref{slowGreen}%
c}\label{slowGreenc}%
\end{align}
With this notation at hand the half-space Green's function can be written in
the form%
\begin{equation}
G_{\operatorname*{half}}:=\operatorname*{e}\nolimits^{-s\tau
_{\operatorname*{illu}}}\Theta_{\nu,s}^{\operatorname*{illu}}%
+\operatorname*{e}\nolimits^{-s\tau_{\operatorname*{refl}}}\left(  \Theta
_{\nu,s}^{\operatorname*{refl}}+\Theta_{\nu,s}^{\operatorname*{imp}}\right)  .
\label{DefGhalf}%
\end{equation}
The functions in (\ref{slowGreen}) are collected in the families\footnote{Note
that the domain of the function $\Theta_{\nu,s}^{\operatorname*{illu}}$ in
(\ref{DefGhalf}) is $H_{+}\times H_{+}$. However, for the analysis of
$\Theta_{\nu,s}^{\operatorname*{refl}}$ it turns out to be useful to consider
$\Theta_{\nu}^{\operatorname*{illu}}$ on the larger domain $\mathbb{R}%
^{d}\times\mathbb{R}^{d}$.}%
\begin{subequations}
\label{slowfamilies}
\end{subequations}%
\begin{align}
\mathcal{F}_{\nu}^{\operatorname*{illu}}  &  :=\left\{  \Theta_{\nu
,s}^{\operatorname*{illu}}:\mathbb{R}^{d}\times\mathbb{R}^{d}\rightarrow
\mathbb{C}\mid s\in\overset{\bullet}{\mathbb{C}}_{\geq0}\right\}  ,\tag{%
\ref{slowfamilies}%
a}\label{slowfamiliesa}\\
\mathcal{F}_{\nu}^{\operatorname*{refl}}  &  :=\left\{  \Theta_{\nu
,s}^{\operatorname*{refl}}:H_{+}\times H_{+}\rightarrow\mathbb{C}\mid
s\in\overset{\bullet}{\mathbb{C}}_{\geq0}\right\}  ,\tag{%
\ref{slowfamilies}%
b}\label{slowfamiliesb}\\
\mathcal{F}_{\nu}^{\operatorname*{imp}}  &  :=\left\{  \Theta_{\nu
,s}^{\operatorname*{imp}}:H_{+}\times H_{+}\rightarrow\mathbb{C}\mid
s\in\overset{\bullet}{\mathbb{C}}_{\geq0}\right\}  \tag{%
\ref{slowfamilies}%
c}\label{slowfamiliesc}%
\end{align}
and our goal is to prove that these families are slowly varying.

\subsection{Families of slowly varying functions\label{SecSlowVaryFct}}

In this section we introduce the notion of families of slowly varying
functions via polynomial approximability, which is the key mechanism of most
numerical approximation methods. Loosely speaking a family of functions
depending on a (frequency) parameter, say $s$, is slowly varying if a
polynomial approximation converges exponentially. Polynomial approximability
on a real interval for \textit{analytic functions} relies on the modulus of
these functions on \textit{Bernstein ellipses} encircling the interval. We
start this section with some elementary geometric considerations and recall
the tensor Chebyshev interpolation.

For $\mathbf{a}=\left(  a_{i}\right)  _{i=1}^{d}\in\mathbb{R}^{d}$ and
$\mathbf{b}=\left(  b_{i}\right)  _{i=1}^{d}\in\mathbb{R}^{d}$ with
$-\infty<a_{i}<b_{i}<\infty$ for all $i\in\left\{  1,2,\ldots,d\right\}  $,
the corresponding cuboid is
\[
\left[  \mathbf{a},\mathbf{b}\right]  :=%
\BIGOP{\times}%
_{j=1}^{d}\left[  a_{j},b_{j}\right]  .
\]

We recall the tensor Chebyshev interpolation on a block of cuboids $\left[
\mathbf{a},\mathbf{b}\right]  \times\left[  \mathbf{c},\mathbf{d}\right]  $.
Let $\hat{\xi}_{i,m}$, $0\leq i\leq m$, denote the Chebyshev points in the
unit interval $\left[  -1,1\right]  $ and let $\hat{L}_{i,m}$ be the
corresponding Lagrange polynomials. The tensor version employs the index set%
\[
\iota_{m}:=\left\{  0,1,\ldots,m\right\}  ^{d}.
\]
and is given, for $%
\mbox{\boldmath$ \mu$}%
\in\iota_{m}$, by $\widehat{%
\mbox{\boldmath$ \xi$}%
}_{%
\mbox{\boldmath$ \mu$}%
,m}:=\left(  \hat{\xi}_{\mu_{1},m},\hat{\xi}_{\mu_{2},m},\ldots,\hat{\xi}%
_{\mu_{d},m}\right)  ^{\intercal}$ and $\hat{L}_{%
\mbox{\boldmath$ \mu$}%
,m}=%
{\displaystyle\bigotimes\limits_{\ell=1}^{d}}
\hat{L}_{\mu_{\ell},m}$. For a box $\left[  \mathbf{a},\mathbf{b}\right]  $,
let $\chi_{\left[  \mathbf{a},\mathbf{b}\right]  }$ denote an affine pullback
to the reference cuboid $\left[  -1,1\right]  ^{d}$. Then, the tensorized
Chebyshev nodal points of order $m$ scaled to the box $\left[  \mathbf{a}%
,\mathbf{b}\right]  $ are given by $%
\mbox{\boldmath$ \xi$}%
_{%
\mbox{\boldmath$ \mu$}%
,m}^{\left[  \mathbf{a},\mathbf{b}\right]  }:=\chi_{\left[  \mathbf{a}%
,\mathbf{b}\right]  }\left(  \widehat{%
\mbox{\boldmath$ \xi$}%
}_{%
\mbox{\boldmath$ \mu$}%
,m}\right)  $ and $L_{%
\mbox{\boldmath$ \mu$}%
}^{\left[  \mathbf{a},\mathbf{b}\right]  }:=\hat{L}_{%
\mbox{\boldmath$ \mu$}%
,m}\circ\chi_{\left[  \mathbf{a},\mathbf{b}\right]  }^{-1}$, for all $%
\mbox{\boldmath$ \mu$}%
\in\mathbb{\iota}_{m}$.

\begin{definition}
\label{DefChebInt}Let $\left[  \mathbf{a},\mathbf{b}\right]  ,\left[
\mathbf{c},\mathbf{d}\right]  $ be two axes-parallel cuboids. The \emph{tensor
Chebyshev interpolation operator} $\Pi_{m}^{\left[  \mathbf{a},\mathbf{b}%
\right]  \times\left[  \mathbf{c} \mathbf{d}\right]  } : C^{0}\left(  \left[
\mathbf{a},\mathbf{b}\right]  \times\left[  \mathbf{c},\mathbf{d}\right]
\right)  \rightarrow{\mathcal{Q}}_{m}:= \operatorname{span}
\{L_{\boldsymbol{\mu},m} \otimes L_{\boldsymbol{\nu},m}\,|\, {\boldsymbol{\mu
}}, {\boldsymbol{\nu}} \in\{0,\ldots,m\}^{d}\}$ of degree $m\in\mathbb{N}_{0}$
is given by%
\[
\Pi_{m}^{\left[  \mathbf{a},\mathbf{b}\right]  \times\left[  \mathbf{c}%
,\mathbf{d}\right]  }k:=\sum_{%
\mbox{\boldmath$ \mu$}%
,%
\mbox{\boldmath$ \nu$}%
\in\iota_{m}}k\left(
\mbox{\boldmath$ \xi$}%
_{%
\mbox{\boldmath$ \mu$}%
,m}^{\left[  \mathbf{a},\mathbf{b}\right]  },%
\mbox{\boldmath$ \xi$}%
_{%
\mbox{\boldmath$ \nu$}%
,m}^{\left[  \mathbf{c},\mathbf{d}\right]  }\right)  L_{%
\mbox{\boldmath$ \mu$}%
,m}^{\left[  \mathbf{a},\mathbf{b}\right]  }\otimes L_{%
\mbox{\boldmath$ \nu$}%
,m}^{\left[  \mathbf{c},\mathbf{d}\right]  }.
\]

\end{definition}

For the analysis of the approximation error we will employ classical error
estimates for Chebyshev interpolation of analytic functions (see \cite{DR}).
This is done by estimating the modulus of analytic functions on
\textit{Bernstein ellipses}, and we recall their basic properties.

Let $-\infty<a<b<\infty$ and consider the real interval $\left[  a,b\right]
$. Let $\mathcal{E}_{a,b}^{\rho}\subset\mathbb{C}$ be the closed ellipse with
focal points $a$, $b$ and semimajor/semiminor axes $\overline{a}$,
$\overline{b}$ given by%
\begin{equation}
\bar{a}=\frac{\rho^{2}+\left(  \frac{b-a}{2}\right)  ^{2}}{2\rho}\geq
\frac{b-a}{2},\qquad\bar{b}=\frac{\rho^{2}-\left(  \frac{b-a}{2}\right)  ^{2}%
}{2\rho}\geq0, \label{formulaabarbbar}%
\end{equation}
where the estimates become an equality if and only the ellipse collapses to
the interval: $\left[  a,b\right]  =\mathcal{E}_{a,b}^{\rho}$ for
$\rho=\left(  b-a\right)  /2$. The sum of the half-axes is given by $\rho
=\bar{a}+\bar{b}$.

For $j\in\left\{  1,2,\ldots,d\right\}  $, the ellipses $\mathcal{E}%
_{a_{j},b_{j}}^{\rho_{j}}$ refer to the coordinate intervals $\left[
a_{j},b_{j}\right]  $ and the semi-axes sums form the vector $%
\mbox{\boldmath$ \rho$}%
:=\left(  \rho_{i}\right)  _{i=1}^{d}$. For $j\in\left\{  1,2,\ldots
,d\right\}  $, we set%
\begin{equation}
\overrightarrow{\mathcal{E}}_{\mathbf{a},\mathbf{b}}^{j}\left(
\mbox{\boldmath$ \rho$}%
\right)  :=I_{1}\times I_{2}\times\cdots\times I_{j-1}\times\mathcal{E}%
_{a_{j},b_{j}}^{\rho_{j}}\times I_{j+1}\times\cdots\times I_{d}.
\label{singleextension}%
\end{equation}
and denote their union by%
\begin{equation}
\overrightarrow{\mathcal{E}}_{\mathbf{a},\mathbf{b},\mathbf{c},\mathbf{d}%
}\left(
\mbox{\boldmath$ \rho$}%
_{1},%
\mbox{\boldmath$ \rho$}%
_{2}\right)  :=%
{\displaystyle\bigcup\limits_{j=1}^{d}}
\left(  \left(  \overrightarrow{\mathcal{E}}_{\mathbf{a},\mathbf{b}}%
^{j}\left(
\mbox{\boldmath$ \rho$}%
_{1}\right)  \times\left[  \mathbf{c},\mathbf{d}\right]  \right)  \cup\left(
\left[  \mathbf{a},\mathbf{b}\right]  \times\overrightarrow{\mathcal{E}%
}_{\mathbf{c},\mathbf{d}}^{j}\left(
\mbox{\boldmath$ \rho$}%
_{2}\right)  \right)  \right)  . \label{defellarrow}%
\end{equation}

The proof of the following proposition can be found in \cite[Thm.~{7.3.6}%
]{SauterSchwab2010}.

\begin{proposition}
[{{\cite[Thm.~{7.3.6}]{SauterSchwab2010}}}]Let $\left[  \mathbf{a}%
,\mathbf{b}\right]  $, $\left[  \mathbf{c},\mathbf{d}\right]  $ be
axes-parallel cuboids. Assume that the function $k\in C^{0}\left(  \left[
\mathbf{a},\mathbf{b}\right]  \times\left[  \mathbf{c},\mathbf{d}\right]
\right)  $ can be extended analytically to $\overrightarrow{\mathcal{E}%
}_{\mathbf{a},\mathbf{b},\mathbf{c},\mathbf{d}}\left(
\mbox{\boldmath$ \rho$}%
_{1},%
\mbox{\boldmath$ \rho$}%
_{2}\right)  $ with $\left(
\mbox{\boldmath$ \rho$}%
_{1}\right)  _{i}>\left(  b_{i}-a_{i}\right)  /2$ and $\left(
\mbox{\boldmath$ \rho$}%
_{2}\right)  _{i}>\left(  d_{i}-c_{i}\right)  /2$, $1\leq i\leq d$ (and is
denoted again by $k$). Then, the Chebyshev interpolant in
Def.~\ref{DefChebInt} satisfies the error estimate%
\[
\left\Vert k-\Pi_{m}^{\left[  \mathbf{a},\mathbf{b}\right]  \times\left[
\mathbf{c},\mathbf{d}\right]  }k\right\Vert _{C^{0}\left(  \left[
\mathbf{a},\mathbf{b}\right]  \times\left[  \mathbf{c},\mathbf{d}\right]
\right)  }\leq C_{\gamma}\gamma^{-m}M_{%
\mbox{\boldmath$ \rho$}%
_{1},%
\mbox{\boldmath$ \rho$}%
_{2}}\left(  k\right)
\]
with the relative extension parameter%
\[
\gamma:=\min\left\{  \min\left\{  \frac{2\left(
\mbox{\boldmath$ \rho$}%
_{1}\right)  _{j}}{b_{j}-a_{j}}:1\leq j\leq d\right\}  ,\min\left\{
\frac{2\left(
\mbox{\boldmath$ \rho$}%
_{2}\right)  _{j}}{d_{j}-c_{j}}:1\leq j\leq d\right\}  \right\}
\]
and%
\[
M_{%
\mbox{\boldmath$ \rho$}%
_{1},%
\mbox{\boldmath$ \rho$}%
_{2}}\left(  k\right)  :=\max_{\left(  \mathbf{x},\mathbf{y}\right)
\in\overrightarrow{\mathcal{E}}_{\mathbf{a},\mathbf{b},\mathbf{c},\mathbf{d}%
}\left(
\mbox{\boldmath$ \rho$}%
_{1},%
\mbox{\boldmath$ \rho$}%
_{2}\right)  }\left\vert k\left(  \mathbf{x},\mathbf{y}\right)  \right\vert .
\]
The constant $C_{\gamma}$ is given by%
\[
C_{\gamma}:=\sqrt{d}2^{d+3/2}\left(  1-\gamma^{-2}\right)  ^{-d}.
\]

\end{proposition}

In our application we consider (Green's) functions that depend on a frequency
parameter $s$. We define the notion \textquotedblleft$\kappa-$slowly
varying\textquotedblright\ for such families of functions and first introduce
an \textit{admissibility condition}.

\begin{definition}
\label{etaadm}Let\footnote{For practical applications, $\eta=\eta_{0}%
\in\left\{  1,2\right\}  $ are usual choices.} $0<\eta_{0}=O\left(  1\right)
$. For $\eta\in\left]  0,\eta_{0}\right]  $, two subsets $B$,$C\subset
\mathbb{R}^{d}$ are $\eta-$\emph{admissible} if%
\begin{equation}
\max\left\{  \operatorname*{diam}B,\operatorname*{diam}C\right\}  \leq
\eta\operatorname{dist}\left(  B,C\right)  . \label{etaadm_cond}%
\end{equation}

\end{definition}

In view of (\ref{DefGhalf}) the first prefactor is oscillatory if
$\operatorname{Im}s\neq0$ and the product \textquotedblleft$\left(
\operatorname{Im}s\right)  \times\left\Vert \mathbf{x}-\mathbf{y}\right\Vert
$\textquotedblright\ becomes large. More specifically, we say that the
parameters $s$, $\mathbf{x},\mathbf{\ y}$ belong to the \textit{slowly
oscillatory regime} if $\left\vert s\right\vert \left\Vert \mathbf{x}%
-\mathbf{y}\right\Vert \leq1$ while they belong to the oscillatory regime as
$\left\vert s\right\vert \left\Vert \mathbf{x}-\mathbf{y}\right\Vert \geq1$
becomes large. These different ranges are reflected in Definition
\ref{DefSlowVar} of a $\kappa-$slowly varying family of functions. First, we
introduce the notation of algebraically bounded functions.

\begin{definition}
Let $\omega\subset\mathbb{C}$. A function $g:\omega\rightarrow\mathbb{R}%
_{\geq0}$ is \emph{algebraically bounded} for growth parameters $%
\mbox{\boldmath$ \alpha$}%
=\left(  \alpha_{1},\alpha_{2}\right)  \in\mathbb{R}^{2}$ and some $C\geq0$ if%
\[
g\left(  \zeta\right)  \leq C\times\left\{
\begin{array}
[c]{cc}%
\left\vert \zeta\right\vert ^{-\alpha_{1}} & \text{if }\left\vert
\zeta\right\vert \geq1,\\
\left\vert \zeta\right\vert ^{-\alpha_{2}} & \text{if }\left\vert
\zeta\right\vert \leq1,
\end{array}
\right.  \quad\forall\zeta\in\omega.
\]

\end{definition}

The definition of a $\kappa-$slowly varying family of functions relies on the
holomorphic extensibility of functions from coordinate intervals to complex
ellipses. While the parameter $\eta$ (cf.\ (\ref{etaadm_cond})) is related to
the admissibility condition for pairs of cuboids, a further parameter
$\kappa>0$ measures the relative size of the extended region: We define for
$j\in\left\{  1,2,\ldots,d\right\}  :$%
\begin{equation}
\left(
\mbox{\boldmath$ \rho$}%
_{1}\right)  _{j}:=\left(
\mbox{\boldmath$ \rho$}%
_{1}\right)  _{j}\left(  \kappa\right)  :=\frac{b_{j}-a_{j}}{2}\left(
1+\frac{2\kappa}{\eta}\right)  ,\text{\quad}\left(
\mbox{\boldmath$ \rho$}%
_{2}\right)  _{j}:=\left(
\mbox{\boldmath$ \rho$}%
_{2}\right)  _{j}\left(  \kappa\right)  :=\frac{d_{j}-c_{j}}{2}\left(
1+\frac{2\kappa}{\eta}\right)  \label{defboldrho}%
\end{equation}
and use the shorthands
\[
\mathcal{E}_{a_{j},b_{j}}\left(  \kappa\right)  :=\mathcal{E}_{a_{j},b_{j}%
}^{({\boldsymbol{\rho}}_{1})_{j}\left(  \kappa\right)  } ,\quad\overrightarrow
{\mathcal{E}}_{\mathbf{a},\mathbf{b}}^{j}\left(  \kappa\right)
:=\overrightarrow{\mathcal{E}}_{\mathbf{a},\mathbf{b}}^{j}\left(
\mbox{\boldmath$ \rho$}%
_{1}\left(  \kappa\right)  \right)  ,\quad\overrightarrow{\mathcal{E}%
}_{\mathbf{a},\mathbf{b},\mathbf{c},\mathbf{d}}\left(  \kappa\right)
:=\overrightarrow{\mathcal{E}}_{\mathbf{a},\mathbf{b},\mathbf{c},\mathbf{d}%
}\left(
\mbox{\boldmath$ \rho$}%
_{1}\left(  \kappa\right)  ,%
\mbox{\boldmath$ \rho$}%
_{2}\left(  \kappa\right)  \right)  .
\]

\begin{definition}
\label{DefSlowVar}Let $0<\eta\leq\eta_{0}$, $\kappa>0$, and $\omega
\subset\mathbb{C}$ be given. For fixed $D\subset\mathbb{R}^{d}$, consider a
family of functions%
\[
\mathcal{F}:=\left\{  F_{s}:D\times D\rightarrow\mathbb{C}\mid s\in
\omega\right\}
\]
along with an algebraically bounded \emph{reference function} $\lambda
:\mathbb{C}\rightarrow\mathbb{R}_{>0}$ and constants $C_{\mathcal{F}}>0$,
$\tau\in\mathbb{R}$.

The family $\mathcal{F}$ is $\kappa-$\emph{slowly varying} if for any block of
$\eta$-admissible cuboids $B:=\left[  \mathbf{a},\mathbf{b}\right]
\times\left[  \mathbf{c},\mathbf{d}\right]  \subset D\times D$ with distance
$\delta:=\operatorname{dist}\left(  \left[  \mathbf{a},\mathbf{b}\right]
,\left[  \mathbf{c},\mathbf{d}\right]  \right)  >0$ and any $F_{s}%
\in{\mathcal{F}}$, the function $\left.  F_{s}\right\vert _{B}$ can be
extended analytically to $\overrightarrow{\mathcal{E}}_{\mathbf{a}%
,\mathbf{b},\mathbf{c},\mathbf{d}}\left(  \kappa\right)  $ and satisfies%
\[
\max_{\left(  \mathbf{x},\mathbf{y}\right)  \in\overrightarrow{\mathcal{E}%
}_{\mathbf{a},\mathbf{b},\mathbf{c},\mathbf{d}}\left(  \kappa\right)
}\left\vert F_{s}\left(  \mathbf{x},\mathbf{y}\right)  \right\vert \leq
C_{\mathcal{F}}\left\vert s\right\vert ^{\tau}\lambda\left(  \left\vert
s\right\vert \delta\right)  .
\]

\end{definition}

For $\kappa-$slowly varying function families the tensor Chebyshev
interpolation converges exponentially as can be seen from the following corollary.

\begin{corollary}
\label{CorMaj}Let $\eta>0$, $\kappa>0,$ and $\omega\subset\mathbb{C}$ be
given. For fixed $D\subseteq\mathbb{R}^{d}$, consider a family of functions
\[
\mathcal{F}:=\left\{  k_{s}:D\times D\rightarrow\mathbb{C}\mid s\in
\omega\right\}
\]
that is $\kappa-$\emph{slowly varying} with reference function $\lambda$ and
constants $C_{\mathcal{F}}$, $\tau$ as in Definition~\ref{DefSlowVar}.

Then, for any $s\in\omega$ and any block of $\eta$-admissible cuboids $\left[
\mathbf{a},\mathbf{b}\right]  \times\left[  \mathbf{c},\mathbf{d}\right]
\subset D\times D$ it holds:%
\[
\left\Vert k_{s}-\Pi_{m}^{\left[  \mathbf{a},\mathbf{b}\right]  \times\left[
\mathbf{c},\mathbf{d}\right]  }k_{s}\right\Vert _{C^{0}\left(  \left[
\mathbf{a},\mathbf{b}\right]  \times\left[  \mathbf{c},\mathbf{d}\right]
\right)  }\leq C_{\mathcal{F}}\left\vert s\right\vert ^{\tau}\lambda\left(
\left\vert s\right\vert \delta\right)  C_{\gamma}\gamma^{-m}\quad\forall
m\in\mathbb{N}_{0}%
\]
with $\gamma=1+2\kappa/\eta$ and $\delta:=\operatorname{dist}\left(  \left[
\mathbf{a},\mathbf{b}\right]  ,\left[  \mathbf{c},\mathbf{d}\right]  \right)
$.

Since the function $\lambda$ is algebraically bounded, a relative accuracy
$\left\Vert k_{s}-\Pi_{m}^{\left[  \mathbf{a},\mathbf{b}\right]  \times\left[
\mathbf{c},\mathbf{d}\right]  }k_{s}\right\Vert _{C^{0}\left(  \left[
\mathbf{a},\mathbf{b}\right]  \times\left[  \mathbf{c},\mathbf{d}\right]
\right)  }\leq\varepsilon\left\vert s\right\vert ^{\tau}\lambda\left(
\left\vert s\right\vert \delta\right)  $ for given $\varepsilon>0$ is reached
for a polynomial degree $m$ which depends only linearly on $\left\vert
\ln\varepsilon\right\vert $, $\left\vert \ln\left\vert s\right\vert
\right\vert $, and $\left\vert \ln\left\vert s\delta\right\vert \right\vert $.
\end{corollary}

\begin{lemma}
\label{lemma:slowly-varying-analyticity}Consider a family of function
$\mathcal{F}:=\left\{  k_{s}:D\times D\rightarrow\mathbb{C}\mid s\in
\omega\right\}  $ that is $\kappa-$\emph{slowly varying} with reference
function $\lambda$ as in Definition~\ref{DefSlowVar} Then for every $\eta
$-admissible cuboid ${\mathbf{B}}=[\mathbf{a},\mathbf{b}]\times\lbrack
\mathbf{c},\mathbf{d}]\subset D\times D$ there is a complex neighborhood
${\mathcal{B}}\subset{\mathbb{C}}^{d}$ of ${\mathbf{B}}$ such that every
$F_{s}\in{\mathcal{F}}$ can be extended analytically to ${\mathcal{B}}$.
Furthermore, given $T>0$ there are constants $C_{1}$, $C_{2}$ depending solely
on $T$ such that for $\kappa/\eta\leq T$ there holds upon setting
$\Theta:=\min_{j=1,\ldots,d}\min\{(b_{j}-a_{j}),(d_{j}-c_{j})\}$
\begin{equation}
\Vert\partial_{x}^{\boldsymbol{\mu}}\partial_{y}^{\boldsymbol{\nu}}F_{s}%
\Vert_{C^{0}(\mathbf{B})}\leq C_{1}C_{\mathcal{F}}|s|^{\tau}\lambda
(s\delta){\boldsymbol{\mu}}!{\boldsymbol{\nu}}!\left(  C_{2}\left(
\eta/\kappa\right)  ^{2}\Theta\right)  ^{|{\boldsymbol{\mu}}%
|+|{\boldsymbol{\nu}}|}. \label{eq:lemma:slowly-varying-analyticity-10}%
\end{equation}

\end{lemma}

\proof
By assumption, the function $F_{s}$ can be extended analytically in each
variable. For quantitative bounds, we have to consider $\operatorname{dist}%
([a_{j},b_{j}], \partial\mathcal{E}_{a_{j},b_{j}}(\kappa))$. For small
$\kappa/\eta$, one has $\operatorname{dist}([a_{j},b_{j}], \partial
\mathcal{E}_{a_{j},b_{j}}(\kappa)) \gtrsim(b_{j}-a_{j}) \left(  \frac{\kappa
}{\eta}\right)  ^{2}$ with implied constant independent of $a_{j}$, $b_{j}$,
$\kappa/\eta$. An analogous result holds for the intervals $[c_{j},d_{j}]$.
Hence, Cauchy's integral theorem allows us to control the partial derivatives
of $F_{s}$ in the stated fashion.
\endproof

\section{Analysis of the families $\Theta_{\nu,s}^{\operatorname*{illu}%
},\Theta_{\nu,s}^{\operatorname*{refl}},\Theta_{\nu,s}^{\operatorname*{imp}}$}

In this chapter we will prove that the families of functions $\Theta_{\nu
,s}^{\operatorname*{illu}}$, $\Theta_{\nu,s}^{\operatorname*{refl}}$,
$\Theta_{\nu,s}^{\operatorname*{imp}}$ in (\ref{slowGreen}) are $\kappa
-$slowly varying. First, we present a majorant for the modified Bessel
function $K_{\nu}$ which is uniform for all $z\in\mathbb{C}$ (see
\S \ref{SecMacdonald1}). Since the functions $\Theta_{\nu,s}%
^{\operatorname*{illu}}$ and $\Theta_{\nu,s}^{\operatorname*{refl}}$ depend on
the Euclidean norm $\left\Vert \mathbf{x}-\mathbf{y}\right\Vert $ and
$\left\Vert \mathbf{x}-\mathbf{Ry}\right\Vert $ we will derive estimates for
the holomorphic norm extension to complex ellipses (see \S \ref{SecHoloNorm}).
These results are then combined to prove that $\Theta_{\nu,s}%
^{\operatorname*{illu}}$ and $\Theta_{\nu,s}^{\operatorname*{refl}}$ are
$\kappa-$slowly varying (\S \ref{SecAnaGreenI}). Finally, in Section
\ref{SecAnaGreenII} we show that $\Theta_{\nu,s}^{\operatorname*{imp}}$ is
$\kappa-$slowly varying and the main theorem (Thm.~\ref{TheoSlowVar}) of this
paper follows.

\begin{notation}
$C_{\nu}$ is a constant depending only on $\nu\in\mathbb{Z}%
^{\operatorname*{half}}$ and may change its value in each appearance.

$C_{\nu,\beta}$ depends on $\nu\in\mathbb{Z}^{\operatorname*{half}}$ and is a
continuous function of the impedance parameter $\beta>0$. It may exhibit
algebraic or logarithmic singularities towards the endpoints $\beta\in\left\{
0,\infty\right\}  $ whose strength depends on $\nu$. It may change its value
in each appearance.
\end{notation}

For the formulation of the main theorem, two (similar) functions $M_{\mu
}:\mathbb{R}_{>0}\rightarrow\mathbb{R}_{>0}$ and $W_{\mu}:\mathbb{R}%
_{>0}\rightarrow\mathbb{R}_{>0}$ are needed%
\begin{subequations}
\label{DefMW}
\end{subequations}%
\begin{align}
M_{\mu}\left(  r\right)   &  :=\left\{
\begin{array}
[c]{ll}%
r^{-1/2} & \text{for }r\geq1,\\
r^{-\mu} & \text{for }0<r\leq1\wedge\mu\in\mathbb{Z}_{\geq\frac{1}{2}%
}^{\operatorname*{half}},\\
1+\left\vert \ln r\right\vert  & \text{for }0<r\leq1\wedge\mu=0,
\end{array}
\right. \tag{%
\ref{DefMW}%
a}\label{majorantdef}\\
W_{\mu}\left(  r\right)   &  :=\left\{
\begin{array}
[c]{ll}%
r^{-1/2} & \text{for }r\geq1,\\
r^{-\mu} & \text{for }0<r\leq1\wedge\mu\in\mathbb{Z}_{\geq\frac{1}{2}%
}^{\operatorname*{half}},\\
1+\ln^{2}r & \text{for }0<r\leq1\wedge\mu=0.
\end{array}
\right.  \tag{%
\ref{DefMW}%
b}\label{DefMWb}%
\end{align}

\begin{theorem}
\label{TheoSlowVar}Let $0<\eta\leq\eta_{0}$ be as in Definition \ref{etaadm},
$\nu\in\mathbb{Z}^{\operatorname*{half}},$ and let $\beta>0$ denote the
impedance parameter in (\ref{goveq}). The families of functions $\mathcal{F}%
_{\nu}^{\operatorname*{illu}}$, $\mathcal{F}_{\nu}^{\operatorname*{refl}}$,
$\mathcal{F}_{\nu}^{\operatorname*{imp}}$ in (\ref{slowfamilies}) are
$\kappa-$slowly varying:

\begin{enumerate}
\item for any $0<\kappa<1/6$, any $\eta$-admissible block $B=\left[
\mathbf{a},\mathbf{b}\right]  \times\left[  \mathbf{c},\mathbf{d}\right]
\subset H_{+}\times H_{+}$ and $\delta:=\operatorname{dist}\left(  \left[
\mathbf{a},\mathbf{b}\right]  ,\left[  \mathbf{c},\mathbf{d}\right]  \right)
$ it holds%
\begin{subequations}
\label{TheoIlluRefl}
\end{subequations}%
\begin{align}
\max_{\left(  \mathbf{x},\mathbf{y}\right)  \in\overrightarrow{\mathcal{E}%
}_{\mathbf{a},\mathbf{b},\mathbf{c},\mathbf{d}}\left(  \kappa\right)
}\left\vert \Theta_{\nu,s}^{\operatorname*{illu}}\left(  \mathbf{x}%
,\mathbf{y}\right)  \right\vert  &  \leq C_{\nu}\left(  \frac{\left\vert
s\right\vert }{\delta}\right)  ^{\nu+1/2}M_{\nu+1/2}\left(  \left\vert
s\right\vert \delta\right)  ,\tag{%
\ref{TheoIlluRefl}%
a}\label{TheoIlluRefla}\\
\max_{\left(  \mathbf{x},\mathbf{y}\right)  \in\overrightarrow{\mathcal{E}%
}_{\mathbf{a},\mathbf{b},\mathbf{c},\mathbf{d}}\left(  \kappa\right)
}\left\vert \Theta_{\nu,s}^{\operatorname*{refl}}\left(  \mathbf{x}%
,\mathbf{y}\right)  \right\vert  &  \leq C_{\nu,\beta}\left(  \frac{\left\vert
s\right\vert }{\delta}\right)  ^{\nu+1/2}M_{\nu+1/2}\left(  \left\vert
s\right\vert \delta\right)  , \tag{%
\ref{TheoIlluRefl}%
b}\label{TheoIlluReflb}%
\end{align}
with $M_{\nu+1/2}$ as in (\ref{majorantdef}).

The constant $\mathfrak{C}_{s}$ in Definition \ref{DefSlowVar} for the first
case (\ref{TheoIlluRefla}) can be chosen as $\mathfrak{C}_{s}:=C_{\nu
}\left\vert s\right\vert ^{2\nu+1}$ and in the second one (\ref{TheoIlluReflb}%
) by $\mathfrak{C}_{s}:=C_{\nu,\beta}\left\vert s\right\vert ^{2\nu+1}$. The
reference functions $\lambda$ can be chosen for both cases by%
\begin{equation}
\lambda:\overset{\bullet}{\mathbb{C}}_{\geq0}\rightarrow\mathbb{R}_{>0}%
\quad\lambda\left(  \zeta\right)  :=\left\{
\begin{array}
[c]{ll}%
\left\vert \zeta\right\vert ^{-\nu-1} & \text{for }\left\vert \zeta\right\vert
\geq1,\\
\left\vert \zeta\right\vert ^{-2\nu-1} & \text{for }0<\left\vert
\zeta\right\vert \leq1\wedge\nu\in\mathbb{Z}_{\geq0}^{\operatorname*{half}},\\
1+\left\vert \ln\left\vert \zeta\right\vert \right\vert  & \text{for
}0<\left\vert \zeta\right\vert \leq1\wedge\nu=-1/2.
\end{array}
\right.  \label{defCslambda_illu}%
\end{equation}

\item There exists a positive number $C_{\mathcal{E}}$ independent of all
parameters and functions such that for any $\kappa\in\left[  0,\frac{\beta
^{2}}{C_{\mathcal{E}}\left(  1+\beta\right)  ^{3}}\right[  $, any $\eta
$-admissible block $B=\left[  \mathbf{a},\mathbf{b}\right]  \times\left[
\mathbf{c},\mathbf{d}\right]  \subset H_{+}\times H_{+}$ with $\delta
:=\operatorname{dist}\left(  \left[  \mathbf{a},\mathbf{b}\right]  ,\left[
\mathbf{c},\mathbf{d}\right]  \right)  $, it holds%
\[
\max_{\left(  \mathbf{x},\mathbf{y}\right)  \in\overrightarrow{\mathcal{E}%
}_{\mathbf{a},\mathbf{b},\mathbf{c},\mathbf{d}}\left(  \kappa\right)
}\left\vert \Theta_{\nu,s}^{\operatorname*{imp}}\left(  \mathbf{x}%
,\mathbf{y}\right)  \right\vert \leq C_{\nu,\beta}\left(  \frac{\left\vert
s\right\vert }{\delta}\right)  ^{\nu+1/2}W_{\nu+1/2}\left(  \left\vert
s\right\vert \delta\right)
\]
with $W_{\mu}$ as in (\ref{DefMWb}). The constant $\mathfrak{C}_{s}$ in
Definition \ref{DefSlowVar} can be chosen by $\mathfrak{C}_{s}:=C_{\nu,\beta
}\left\vert s\right\vert ^{2\nu+1}$ and the reference function by%
\[
\lambda:\overset{\bullet}{\mathbb{C}}_{\geq0}\rightarrow\mathbb{R}_{>0}%
\quad\lambda\left(  \zeta\right)  :=\left\{
\begin{array}
[c]{ll}%
\left\vert \zeta\right\vert ^{-\nu-1} & \text{for }\left\vert \zeta\right\vert
\geq1,\\
\left\vert \zeta\right\vert ^{-2\nu-1} & \text{for }0<\left\vert
\zeta\right\vert \leq1\wedge\nu\in\mathbb{Z}_{\geq0}^{\operatorname*{half}},\\
1+\ln^{2}\left\vert \zeta\right\vert  & \text{for }0<\left\vert \zeta
\right\vert \leq1\wedge\nu=-1/2.
\end{array}
\right.
\]

\end{enumerate}
\end{theorem}

The proof of this theorem is a consequence of Lemmata \ref{TheoIlluRefllem}
and \ref{LemVarimpII}.

\subsection{A majorant for the Macdonald function\label{SecMacdonald1}}

In this section, we derive a majorant for the modified Bessel function of
second kind. For $R>0$, let $B_{R}:=\left\{  \zeta\in\mathbb{C}\mid\left\vert
\zeta\right\vert \leq R\right\}  $.

\begin{lemma}
\label{LemKmue}Let $\mu\in\mathbb{Z}_{\geq0}^{\operatorname*{half}}$. The
modified Bessel function $K_{\mu}$ satisfies the estimates%
\begin{subequations}
\label{Kmue_est}
\end{subequations}%

\begin{enumerate}
\item[a.]
\begin{equation}
\left\vert \operatorname*{e}\nolimits^{z}K_{\mu}\left(  z\right)  \right\vert
\leq\sqrt{\frac{\pi}{2\left\vert z\right\vert }}\exp\left(  \frac{\left\vert
\mu^{2}-\tfrac{1}{4}\right\vert }{\left\vert z\right\vert }\right)
\quad\forall z\in\mathbb{C}\backslash\mathbb{R}_{\leq0}, \tag{%
\ref{Kmue_est}%
a}\label{Kmue_esta}%
\end{equation}

\item[b.]
\begin{equation}
\left\vert \operatorname*{e}\nolimits^{z}\left(  K_{\mu}\left(  z\right)
-K_{\mu+1}\left(  z\right)  \right)  \right\vert \leq\frac{C_{\mu}}{\left\vert
z\right\vert ^{3/2}}\exp\left(  \frac{C_{\mu}}{\left\vert z\right\vert
}\right)  \quad\forall z\in\mathbb{C}\backslash\mathbb{R}_{\leq0} \tag{%
\ref{Kmue_est}%
b}\label{Kmue_esta2}%
\end{equation}
with $C_{\mu}:=\sqrt{8\pi^{3}}\left(  \mu+1\right)  ^{2}$.

\item[c.] For any $R>0$, there exists a constant $C_{\mu}\left(  R\right)  $
depending only on $\mu$ and $R$ such that the Bessel function $K_{\mu}$
satisfies:%
\begin{equation}
\left\vert \operatorname*{e}\nolimits^{z}K_{\mu}\left(  z\right)  \right\vert
\leq C_{\mu}\left(  R\right)  \times\left\{
\begin{array}
[c]{ll}%
\left(  1+\left\vert \ln\left\vert z\right\vert \right\vert \right)  &
\mu=0,\\
\left(  \dfrac{\left\vert z\right\vert }{2}\right)  ^{-\mu} &
\text{otherwise,}%
\end{array}
\right.  \quad\forall z\in B_{R}\backslash\left(  \left[  -R,0\right]
\right)  . \tag{%
\ref{Kmue_est}%
c}\label{Kmue_estb}%
\end{equation}

\end{enumerate}
\end{lemma}%

\proof
Estimate (\ref{Kmue_esta}) follows from the well-known asymptotic expansion of
the modified Bessel function for large argument (see, \cite[(86)]%
{Gergoe_K_nue_bound} with reference to \cite[p. 270]{OlverAsymptotics}). For
(\ref{Kmue_esta2}) we use the asymptotic expansion from \cite[10.40.10, .11,
.12]{NIST:DLMF}%

\[
\operatorname*{e}\nolimits^{z}\left(  K_{\mu}\left(  z\right)  -K_{\mu
+1}\left(  z\right)  \right)  =\left(  \frac{\pi}{2z}\right)  ^{\frac{1}{2}%
}\left(  1+R_{1}\left(  \mu,z\right)  -\left(  1+R_{1}\left(  \mu+1,z\right)
\right)  \right)
\]
with\footnote{Here, rather generously, the maximal prefactor $2\chi\left(
1\right)  =\pi$ for the estimate of the variational operator in
\cite[10.40.12]{NIST:DLMF} has been used.}%

\[
\left\vert R_{1}\left(  \nu,z\right)  \right\vert \leq\frac{\pi}{4}\left\vert
\frac{4\nu^{2}-1}{z}\right\vert \exp\left(  \pi\frac{\left\vert \nu^{2}%
-\tfrac{1}{4}\right\vert }{\left\vert z\right\vert }\right)  .
\]
Hence%
\begin{align*}
\left\vert \operatorname*{e}\nolimits^{z}\left(  K_{\mu}\left(  z\right)
-K_{\mu+1}\left(  z\right)  \right)  \right\vert  &  \leq\left(  \frac{\pi
}{2\left\vert z\right\vert }\right)  ^{\frac{1}{2}}\left(  \left\vert
R_{1}\left(  \mu,z\right)  \right\vert +\left\vert R_{1}\left(  \mu
+1,z\right)  \right\vert \right) \\
&  \leq\left(  \frac{2\pi}{\left\vert z\right\vert }\right)  ^{\frac{1}{2}%
}\left\vert R_{1}\left(  \mu+1,z\right)  \right\vert \leq\frac{C_{\mu}%
}{\left\vert z\right\vert ^{3/2}}\exp\left(  \frac{C_{\mu}}{\left\vert
z\right\vert }\right)  .
\end{align*}

Next we prove (\ref{Kmue_estb}) and start with some preparations. Recall the
$n$-th harmonic number%
\[
H_{n}:=\sum_{m=1}^{n}\frac{1}{m}%
\]
and the Euler--Mascheroni constant $\gamma=\lim_{n\rightarrow+\infty}\left(
H_{n}-\ln n\right)  =0.57721\;56649\;01532\;86060\;\dots.$. For $\nu
\in\mathbb{N}$, the combination of \cite[10.31.1]{NIST:DLMF} with
\cite[10.25.2]{NIST:DLMF} leads to%
\begin{subequations}
\label{Knueint}
\end{subequations}%
\begin{equation}
K_{\nu}\left(  z\right)  =\left(  \dfrac{z}{2}\right)  ^{-\nu}\sum
_{k=0}^{\infty}c_{\nu,k}\left(  z\right)  \left(  \dfrac{z}{2}\right)  ^{2k}
\tag{%
\ref{Knueint}%
a}\label{Knueinta}%
\end{equation}
for%
\begin{equation}
c_{\nu,k}\left(  z\right)  :=\left\{
\begin{array}
[c]{ll}%
\dfrac{\left(  -1\right)  ^{k}}{2}\dfrac{(\nu-1-k)!}{k!} & k\in\left\{
0,1,\ldots,\nu-1\right\}  ,\\
\left(  -\ln\dfrac{z}{2}-\gamma+\frac{H_{k-\nu}+H_{k}}{2}\right)
\frac{\left(  -1\right)  ^{\nu}}{\left(  k-\nu\right)  !k!} & k\in\left\{
\nu,\nu+1,\ldots\right\}  .
\end{array}
\right.  \tag{%
\ref{Knueint}%
b}\label{Knueintb}%
\end{equation}
First, we derive bounds for the function $K_{\nu}$ with $\nu\in\mathbb{N}_{0}%
$. From (\ref{Knueint}) we derive%
\begin{equation}
K_{\nu}\left(  z\right)  =\left(  \dfrac{z}{2}\right)  ^{-\nu}\chi_{\nu
}^{\operatorname*{I}}\left(  z\right)  -\left(  \gamma+\ln\dfrac{z}{2}\right)
\left(  \dfrac{z}{2}\right)  ^{\nu}\chi_{\nu}^{\operatorname*{II}}\left(
z\right)  +\left(  \dfrac{z}{2}\right)  ^{\nu}\chi_{\nu}^{\operatorname*{III}%
}\left(  z\right)  \label{Knueanalytic}%
\end{equation}
with the functions%
\begin{align*}
\chi_{\nu}^{\operatorname*{I}}\left(  z\right)   &  :=\sum_{k=0}^{\nu-1}%
\dfrac{\left(  -1\right)  ^{k}}{2}\dfrac{(\nu-1-k)!}{k!}\left(  \dfrac{z}%
{2}\right)  ^{2k},\\
\chi_{\nu}^{\operatorname*{II}}\left(  z\right)   &  :=\sum_{k=0}^{\infty
}\frac{\left(  -1\right)  ^{\nu}}{k!\left(  k+\nu\right)  !}\left(  \dfrac
{z}{2}\right)  ^{2k},\\
\chi_{\nu}^{\operatorname*{III}}\left(  z\right)   &  :=\sum_{k=0}^{\infty
}\frac{H_{k}+H_{k+\nu}}{2}\frac{\left(  -1\right)  ^{\nu}}{k!\left(
k+\nu\right)  !}\left(  \dfrac{z}{2}\right)  ^{2k}.
\end{align*}
All three functions $\chi_{\nu}^{\operatorname*{I},\operatorname*{II}%
,\operatorname*{III}}$ are analytic and satisfy%
\[
\chi_{0}^{\operatorname*{I}}=0\quad\text{and\quad}\max_{z\in B_{R}}%
\max\left\{  \left\vert \chi_{\nu}^{\operatorname*{I}}\left(  z\right)
\right\vert ,\left\vert \chi_{\nu}^{\operatorname*{II}}\left(  z\right)
\right\vert \right\}  \leq\tilde{C}_{\nu}\left(  R\right)
\]
for some constant $0<\tilde{C}_{\nu}\left(  R\right)  <\infty$ which depends
only on $\nu$ and $R$. For $\chi_{\nu}^{\operatorname*{III}}\left(  z\right)
$ we use $\left(  H_{k}+H_{k+\nu}\right)  /2\leq H_{k+\nu}$ to obtain%
\[
\frac{\left(  H_{k}+H_{k+\nu}\right)  /2}{\left(  k+\nu\right)  !}\leq
\frac{H_{k+\nu}}{\left(  k+\nu\right)  !}\leq\frac{1+\ln\left(  1+k+\nu
\right)  }{\left(  k+\nu\right)  !}\leq2.
\]
Hence for $z\in B_{R}$, we get%
\[
\left\vert \chi_{\nu}^{\operatorname*{III}}\left(  z\right)  \right\vert
=\left\vert \sum_{k=0}^{\infty}\frac{H_{k}+H_{k+\nu}}{2}\frac{\left(
-1\right)  ^{\nu}}{k!\left(  k+\nu\right)  !}\left(  \dfrac{z}{2}\right)
^{2k}\right\vert \leq2\left\vert \sum_{k=0}^{\infty}\frac{1}{k!}\left(
\dfrac{R}{2}\right)  ^{2k}\right\vert \leq2\operatorname*{e}\nolimits^{\left(
R/2\right)  ^{2}}=:C_{R}.
\]
For the estimate of $K_{\nu}$ an inequality of the complex logarithm is
needed. We write $z\in\mathbb{C}$ in polar coordinates $z=r\operatorname*{e}%
^{\operatorname*{i}\varphi}$, $r\in\mathbb{R}_{>0}$ and $\varphi\in\left]
-\pi,\pi\right]  $ and get%
\[
\ln z=\ln r+\operatorname*{i}\varphi\quad\text{so that\quad}\left\vert \ln
z\right\vert \leq\left\vert \ln r\right\vert +\pi.
\]
This leads to the estimate of $K_{\nu}:$%
\[
\left\vert K_{\nu}\left(  z\right)  \right\vert \leq\left\vert \dfrac{z}%
{2}\right\vert ^{-\nu}\tilde{C}_{\nu}\left(  R\right)  +\left(  \gamma
+\pi+\left\vert \ln\left\vert \frac{z}{2}\right\vert \right\vert \right)
\left\vert \dfrac{z}{2}\right\vert ^{\nu}\tilde{C}_{\nu}\left(  R\right)
+\left\vert \dfrac{z}{2}\right\vert ^{\nu}C_{R}.
\]
In this way, the function $K_{0}$ can be estimated, for $\nu=0$, by%
\[
\left\vert K_{0}\left(  z\right)  \right\vert \leq C_{0}\left(  R\right)
\left(  1+\left\vert \ln\left\vert \frac{z}{2}\right\vert \right\vert \right)
\quad\forall z\in B_{R}.
\]
For $\nu\in\left\{  1,2,\ldots\right\}  $, we use the fact that powers decay
fast than any logarithmic growth to obtain%
\[
\left(  \gamma+\pi+\left\vert \ln\left\vert \frac{z}{2}\right\vert \right\vert
\right)  \left\vert \dfrac{z}{2}\right\vert ^{\nu}\leq\hat{C}_{\nu}\left(
R\right)  \quad\forall z\in\overset{\bullet}{\mathbb{C}}_{\geq0}\cap B_{R}%
\]
for some constant $\hat{C}_{\nu}\left(  R\right)  $ and obtain%
\begin{equation}
\left\vert K_{\nu}\left(  z\right)  \right\vert \leq C_{\nu}\left(  R\right)
\left(  \dfrac{\left\vert z\right\vert }{2}\right)  ^{-\nu}\quad\forall
z\in\overset{\bullet}{\mathbb{C}}_{\geq0}\cap B_{R}. \label{Knue}%
\end{equation}
For $z\in B_{R}$ the modulus of the prefactor $\operatorname*{e}^{z}$ in
(\ref{Kmue_estb}) is bounded by $\operatorname*{e}^{R}$ so that
(\ref{Kmue_estb}) follows.

Next we prove (\ref{Kmue_estb}) for half integers $\mu=\nu+1/2$, $\nu
\in\left\{  -1,0,1,\ldots\right\}  $. The modified Bessel functions of half
integer have a finite representation (see, e.g., \cite[Chap. 2, (5)]%
{Bessel_polys}):%
\[
\operatorname*{e}\nolimits^{z}K_{\nu+1/2}\left(  z\right)  =\sqrt{\frac{\pi
}{2}}\frac{\theta_{\nu}\left(  z\right)  }{z^{\nu+1/2}}%
\]
with the reverse Bessel polynomials $\theta_{n}$ (see, e.g., \cite[18.34.2]%
{NIST:DLMF}, \cite[Chap. 2, (5)]{Bessel_polys}). Since $\theta_{\nu}$ is a
polynomial of degree $\nu$, estimate (\ref{Kmue_estb}) also holds for half
integers $\mu=\nu+1/2$ and an adjusted constant $C_{\mu}\left(  R\right)  $
still depending only on $\nu$ and $R$.%
\endproof

By considering (\ref{Kmue_esta}), (\ref{Kmue_esta2}) for $\left\vert
z\right\vert \geq1$ and choosing $R=1$ in (\ref{Kmue_estb}) we obtain the
following corollary.

\begin{corollary}
\label{CorFinalMmue}The modified Bessel function $K_{\mu}$ satisfies%
\begin{equation}
\left\vert \operatorname*{e}\nolimits^{z}K_{\mu}\left(  z\right)  \right\vert
\leq C_{\mu}M_{\mu}\left(  \left\vert z\right\vert \right)  \qquad\forall
z\in\mathbb{C}\backslash\mathbb{R}_{\leq0}, \label{majorantest}%
\end{equation}
with $M_{\mu}$ as in (\ref{majorantdef}).

Let $\mu\in\mathbb{Z}_{\geq0}^{\operatorname*{half}}$. Then%
\begin{equation}
\left\vert \operatorname*{e}\nolimits^{z}\left(  K_{\mu}\left(  z\right)
-K_{\mu+1}\left(  z\right)  \right)  \right\vert \leq C_{\mu}\frac{N_{\mu
}\left(  \left\vert z\right\vert \right)  }{\left\vert z\right\vert }%
\quad\forall z\in\mathbb{C}\backslash\mathbb{R}_{\leq0},
\label{mjorantdiffKnue}%
\end{equation}
where the majorant $N_{\mu}:\mathbb{R}_{>0}\rightarrow\mathbb{R}_{>0}$ is
given by%
\[
N_{\mu}\left(  r\right)  :=\left\{
\begin{array}
[c]{ll}%
r^{-1/2} & \text{for }r\geq1,\\
r^{-\mu} & \text{for }0<r\leq1.
\end{array}
\right.
\]

The functions $M_{\mu}$, $N_{\mu}$ are strictly decreasing on $\mathbb{R}%
_{>0}$ and it holds%
\begin{equation}
N_{\mu}\left(  r\right)  \leq M_{\mu}\left(  r\right)  \quad\forall
r\in\mathbb{R}_{>0}. \label{NmueMmue}%
\end{equation}

\end{corollary}%

\proof
The estimates are direct consequences of Lemma \ref{LemKmue} by considering
(\ref{Kmue_esta}), (\ref{Kmue_esta2}) for $\left\vert z\right\vert \geq1$ and
choosing $R=1$ in (\ref{Kmue_estb}). The estimate of $\left\vert
\operatorname*{e}\nolimits^{z}\left(  K_{\mu}\left(  z\right)  -K_{\mu
+1}\left(  z\right)  \right)  \right\vert $ for $\left\vert z\right\vert
\leq1$ follows by a triangle inequality%
\[
\left\vert \operatorname*{e}\nolimits^{z}\left(  K_{\mu}\left(  z\right)
-K_{\mu+1}\left(  z\right)  \right)  \right\vert \leq\left\vert
\operatorname*{e}\nolimits^{z}K_{\mu}\left(  z\right)  \right\vert +\left\vert
\operatorname*{e}\nolimits^{z}K_{\mu+1}\left(  z\right)  \right\vert
\]
in combination with (\ref{Kmue_estb}).%
\endproof

\subsection{Holomorphic norm extension\label{SecHoloNorm}}

The functions in (\ref{slowGreen}) depend on the Euclidean norm $\left\Vert
\mathbf{x}-\mathbf{y}\right\Vert $ and $\left\Vert \mathbf{x}-\mathbf{Ry}%
\right\Vert $; to prove that they are $\kappa-$slowly varying requires to
study holomorphic norm extensions and some preliminaries.

The bilinear form $\left\langle \cdot,\cdot\right\rangle :\mathbb{C}^{d}%
\times\mathbb{C}^{d}\rightarrow\mathbb{C}$ is given, for $\mathbf{u}=\left(
u_{j}\right)  _{j=1}^{d}$, $\mathbf{v}=\left(  v_{j}\right)  _{j=1}^{d}$, by
$\left\langle \mathbf{u},\mathbf{v}\right\rangle =\sum_{j=1}^{d}u_{j}v_{j}$.
For $z\in\overset{\bullet}{\mathbb{C}}_{\geq0}$ we choose $\arg z\in\left[
-\frac{\pi}{2},\frac{\pi}{2}\right]  $ so that $z=\left\vert z\right\vert
\operatorname*{e}^{\operatorname*{i}\arg z}$ and recall the definition of the
principal square root:%
\[
\sqrt{z}:=\sqrt{\left\vert z\right\vert }\operatorname*{e}\nolimits^{\left(
\operatorname*{i}\arg z\right)  /2}.
\]
This directly implies the relations%
\begin{equation}
\left\vert \sqrt{z}\right\vert =\sqrt{\left\vert z\right\vert }\quad
\text{and\quad}\operatorname{Re}\sqrt{z}=\sqrt{\left\vert z\right\vert }%
\cos\left(  \frac{\arg z}{2}\right)  . \label{normz2x}%
\end{equation}
The real part can be estimated from below by\footnote{Here the property
$\arctan x\leq x$ for $x\geq0$ is used.}%
\begin{equation}
\operatorname{Re}\sqrt{z}\geq\sqrt{\left\vert z\right\vert }\left(  1-\frac
{1}{2}\left(  \frac{\arg z}{2}\right)  ^{2}\right)  \geq\sqrt{\left\vert
z\right\vert }\left(  1-\frac{1}{8}\left(  \frac{\operatorname{Im}%
z}{\operatorname{Re}z}\right)  ^{2}\right)  . \label{Reestfrombelow}%
\end{equation}
For the imaginary part we get%
\begin{equation}
\left\vert \operatorname{Im}\sqrt{z}\right\vert =\sqrt{\left\vert z\right\vert
}\left\vert \sin\left(  \frac{\arg z}{2}\right)  \right\vert \leq
\sqrt{\left\vert z\right\vert }\left\vert \sin\left(  \frac{\left\vert
\operatorname{Im}z\right\vert }{2\operatorname{Re}z}\right)  \right\vert
\leq\frac{\left\vert \operatorname{Im}z\right\vert }{2\operatorname{Re}z}%
\sqrt{\left\vert z\right\vert }. \label{estimpart}%
\end{equation}
In the following, we will frequently use the elementary relation
\begin{equation}
z^{2}-\left\vert z\right\vert ^{2}=2\operatorname*{i}z\operatorname{Im}%
z\quad\forall z\in\mathbb{C}\text{.} \label{squarerel}%
\end{equation}

\begin{lemma}
\label{LemExtNorm}Let $\left[  \mathbf{a},\mathbf{b}\right]  \times\left[
\mathbf{c},\mathbf{d}\right]  \subset\mathbb{R}^{d}\times\mathbb{R}^{d}$ be a
block of axes-parallel boxes which are $\eta$-admissible (cf.
(\ref{etaadm_cond})). Let
\[
\kappa\in\left[  0,\kappa_{0}\right[  \quad\text{with\quad}\kappa_{0}:=1/4
\]
and the semi-axes sums be defined by (\ref{defboldrho}). Then, the real part
$\operatorname{Re}\left\langle \mathbf{x}-\mathbf{y},\mathbf{x}-\mathbf{y}%
\right\rangle $ belongs to $\overset{\bullet}{\mathbb{C}}_{\geq0}$ and the
analytic extension of $r$ in (\ref{defn}) to $\overrightarrow{\mathcal{E}%
}_{\mathbf{a},\mathbf{b},\mathbf{c},\mathbf{d}}\left(  \kappa\right)  $, given
by%
\begin{equation}
r\left(  \mathbf{x},\mathbf{y}\right)  :=\sqrt{\left\langle \mathbf{x}%
-\mathbf{y},\mathbf{x}-\mathbf{y}\right\rangle }\quad\forall\left(
\mathbf{x},\mathbf{y}\right)  \in\overrightarrow{\mathcal{E}}_{\mathbf{a}%
,\mathbf{b},\mathbf{c},\mathbf{d}}\left(  \kappa\right)  , \label{defr}%
\end{equation}
is well defined. The extension satisfies the estimates%
\begin{equation}
\left\vert r\left(  \mathbf{x},\mathbf{y}\right)  \right\vert \leq\left\Vert
\mathbf{x}-\mathbf{y}\right\Vert \quad\text{and\quad}\operatorname{Re}r\left(
\mathbf{x},\mathbf{y}\right)  \geq\left(  1-12\kappa^{2}\right)  \left\Vert
\mathbf{x}-\mathbf{y}\right\Vert \quad\forall\left(  \mathbf{x},\mathbf{y}%
\right)  \in\overrightarrow{\mathcal{E}}_{\mathbf{a},\mathbf{b},\mathbf{c}%
,\mathbf{d}}\left(  \kappa\right)  . \label{rests1}%
\end{equation}
For the imaginary part, it holds%
\begin{equation}
\left\vert \operatorname*{Im}r\left(  \mathbf{x},\mathbf{y}\right)
\right\vert \leq4\kappa\left\Vert \mathbf{x}-\mathbf{y}\right\Vert
\quad\forall\left(  \mathbf{x},\mathbf{y}\right)  \in\overrightarrow
{\mathcal{E}}_{\mathbf{a},\mathbf{b},\mathbf{c},\mathbf{d}}\left(
\kappa\right)  . \label{imr1}%
\end{equation}

\end{lemma}

\begin{remark}
In the following estimates we work out the dominant dependence with respect to
$\kappa$; the prefactors are often estimated in a generous way in order to
reduce technicalities in the notation, e.g., for $\kappa\in\left[
0,\kappa_{0}\right[  $ we estimate $\frac{1}{1-\kappa}\leq\frac{4}{3}\leq2$.
However, we emphasize that our estimates are always strict for the full range
of parameters and not only valid \textquotedblleft up to higher order
terms\textquotedblright.
\end{remark}

\textbf{Proof of Lemma \ref{LemExtNorm}. }The upper bound in (\ref{rests1})
follows by a Cauchy-Schwarz inequality and the first relation in
(\ref{normz2x}).

Next we prove the lower bound in (\ref{rests1}). We consider exemplarily the
case that the last side $\left[  a_{d},b_{d}\right]  $ of $\left[
\mathbf{a},\mathbf{b}\right]  $ is extended to $\mathcal{E}_{a_{d},b_{d}%
}\left(  \kappa\right)  $. Let $\left(  \mathbf{x},\mathbf{y}\right)
\in\overrightarrow{\mathcal{E}}_{\mathbf{a},\mathbf{b}}^{d}\left(
\kappa\right)  \times\left[  \mathbf{c},\mathbf{d}\right]  $. Choose $\xi
_{d}\in\left[  a_{d},b_{d}\right]  $ such that%
\begin{equation}
\left\vert x_{d}-\xi_{d}\right\vert =\operatorname{dist}\left(  x_{d},\left[
a_{d},b_{d}\right]  \right)  \quad\text{and set\quad}\mathbf{x}%
_{\operatorname*{R}}:=\left(  \mathbf{x}^{\prime},\xi_{d}\right)  \in\left[
\mathbf{a},\mathbf{b}\right]  . \label{defxid}%
\end{equation}
We will use the shorthands%
\begin{equation}
\mathfrak{D}:=\max\left\{  \operatorname*{diam}\left[  \mathbf{a}%
,\mathbf{b}\right]  ,\operatorname*{diam}\left[  \mathbf{c},\mathbf{d}\right]
\right\}  \quad\text{and\quad}\delta:=\operatorname{dist}\left(  \left[
\mathbf{a},\mathbf{b}\right]  ,\left[  \mathbf{c},\mathbf{d}\right]  \right)
. \label{deltad}%
\end{equation}
and conclude from (\ref{formulaabarbbar}) that%
\begin{align}
\left\vert \xi_{d}-x_{d}\right\vert  &  \leq\max\left\{  \overline{b_{d}%
},\overline{a_{d}}-\frac{b_{d}-a_{d}}{2}\right\}  =\frac{\rho_{d}^{2}-\left(
\frac{b_{d}-a_{d}}{2}\right)  ^{2}}{2\rho_{d}}\label{xidxd}\\
&  \overset{\text{(\ref{defboldrho})}}{=}\frac{\kappa}{\eta}\left(
b_{d}-a_{d}\right)  \frac{\eta+\kappa}{\eta+2\kappa}\leq\frac{\kappa}{\eta
}\left(  b_{d}-a_{d}\right) \nonumber\\
&  \leq\frac{\kappa}{\eta}\mathfrak{D}\leq\kappa\delta.\nonumber
\end{align}
From a triangle inequality it follows%
\begin{equation}
\left\Vert \mathbf{x}-\mathbf{y}\right\Vert \geq\left\Vert \mathbf{x}%
_{\operatorname*{R}}-\mathbf{y}\right\Vert -\left\Vert \mathbf{x}%
_{\operatorname*{R}}-\mathbf{x}\right\Vert \geq\delta-\left\vert \xi_{d}%
-x_{d}\right\vert \label{normfrombelow}%
\end{equation}
and the combination with (\ref{xidxd}) leads to%
\begin{equation}
\left\Vert \mathbf{x}-\mathbf{y}\right\Vert \geq\left(  1-\kappa\right)
\delta. \label{xmylowb}%
\end{equation}
This implies for $\left\vert \xi_{d}-x_{d}\right\vert $ the estimate%
\begin{equation}
\left\vert \xi_{d}-x_{d}\right\vert \leq\frac{\kappa}{1-\kappa}\left\Vert
\mathbf{x}-\mathbf{y}\right\Vert \leq2\kappa\left\Vert \mathbf{x}%
-\mathbf{y}\right\Vert . \label{xidfinal}%
\end{equation}
Next, we investigate the real and imaginary part of $\left\langle
\mathbf{x}-\mathbf{y},\mathbf{x}-\mathbf{y}\right\rangle $. We set
$\omega:=\left\Vert \mathbf{x}^{\prime}-\mathbf{y}^{\prime}\right\Vert $ and
$z_{d}:=x_{d}-y_{d}$ and obtain the relations%
\begin{subequations}
\label{RemImEst0}
\end{subequations}%
\begin{align}
\operatorname{Re}\left\langle \mathbf{x}-\mathbf{y},\mathbf{x}-\mathbf{y}%
\right\rangle  &  =\omega^{2}+\operatorname{Re}\left(  z_{d}^{2}\right)
\overset{\text{(\ref{squarerel})}}{=}\left\Vert \mathbf{x}-\mathbf{y}%
\right\Vert ^{2}-2\left(  \operatorname{Im}z_{d}\right)  ^{2},\tag{%
\ref{RemImEst0}%
a}\label{RemImEst0a}\\
\operatorname{Im}\left\langle \mathbf{x}-\mathbf{y},\mathbf{x}-\mathbf{y}%
\right\rangle  &  \overset{\text{(\ref{squarerel})}}{=}2\operatorname{Re}%
z_{d}\operatorname{Im}z_{d}. \tag{%
\ref{RemImEst0}%
b}\label{RemImEst0b}%
\end{align}
For the imaginary and real part of $z_{d}$, it holds%
\begin{align}
\left\vert \operatorname{Im}z_{d}\right\vert  &  =\left\vert \operatorname{Im}%
\left(  x_{d}-\xi_{d}\right)  \right\vert \leq\left\vert \xi_{d}%
-x_{d}\right\vert \leq2\kappa\left\Vert \mathbf{x}-\mathbf{y}\right\Vert
,\label{estimzd}\\
\left\vert \operatorname{Re}z_{d}\right\vert  &  \leq\left\Vert \mathbf{x}%
-\mathbf{y}\right\Vert , \label{realpartfromabove}%
\end{align}
so that%
\begin{subequations}
\label{ReImEst}
\end{subequations}%
\begin{align}
\operatorname{Re}\left\langle \mathbf{x}-\mathbf{y},\mathbf{x}-\mathbf{y}%
\right\rangle  &  \geq\left(  1-8\kappa^{2}\right)  \left\Vert \mathbf{x}%
-\mathbf{y}\right\Vert ^{2}\overset{\kappa<1/4}{\geq}\frac{\left\Vert
\mathbf{x}-\mathbf{y}\right\Vert ^{2}}{2},\tag{%
\ref{ReImEst}%
a}\label{ReImEsta}\\
\left\vert \operatorname{Im}\left\langle \mathbf{x}-\mathbf{y},\mathbf{x}%
-\mathbf{y}\right\rangle \right\vert  &  \leq4\kappa\left\Vert \mathbf{x}%
-\mathbf{y}\right\Vert ^{2} \tag{%
\ref{ReImEst}%
b}\label{ReImEstb}%
\end{align}
follow. Note that the condition $\kappa\in\left[  0,\kappa_{0}\right[  $
ensures that $\operatorname{Re}\left\langle \mathbf{x}-\mathbf{y}%
,\mathbf{x}-\mathbf{y}\right\rangle >0$ and the function $r$ is well defined
by (\ref{defr}). Hence, we may apply (\ref{ReImEst}) in (\ref{Reestfrombelow})
for%
\begin{align}
\operatorname{Re}\sqrt{\left\langle \mathbf{x}-\mathbf{y},\mathbf{x}%
-\mathbf{y}\right\rangle }  &  \geq\left(  1-8\kappa^{2}\right)
\sqrt{\left\vert \left\langle \mathbf{x}-\mathbf{y},\mathbf{x}-\mathbf{y}%
\right\rangle \right\vert }\nonumber\\
&  \geq\left(  1-8\kappa^{2}\right)  \sqrt{\operatorname{Re}\left\langle
\mathbf{x}-\mathbf{y},\mathbf{x}-\mathbf{y}\right\rangle }\nonumber\\
&  \geq\left(  1-8\kappa^{2}\right)  ^{3/2}\left\Vert \mathbf{x}%
-\mathbf{y}\right\Vert . \label{Refinal}%
\end{align}
Again, $\kappa\in\left[  0,\kappa_{0}\right[  $ implies that the $\kappa
$-dependent factor in the right-hand side of (\ref{Refinal}) is positive. A
straightforward calculus shows that%
\[
\left(  1-8\kappa^{2}\right)  ^{3/2}\geq1-12\kappa^{2}\quad\forall\kappa
\in\left[  0,\kappa_{0}\right[  \text{.}%
\]
For the imaginary part of $r$ we employ (\ref{estimpart}) to get%
\[
\left\vert \operatorname{Im}\sqrt{\left\langle \mathbf{x}-\mathbf{y}%
,\mathbf{x}-\mathbf{y}\right\rangle }\right\vert \overset
{\text{(\ref{ReImEstb})}}{\leq}2\sqrt{\kappa}\left\Vert \mathbf{x}%
-\mathbf{y}\right\Vert \overset{\kappa<1/4}{\leq}4\kappa\left\Vert
\mathbf{x}-\mathbf{y}\right\Vert .
\]%
\endproof

\begin{lemma}
\label{Lemrrplus}Let $\left[  \mathbf{a},\mathbf{b}\right]  \times\left[
\mathbf{c},\mathbf{d}\right]  \subset\mathbb{R}^{d}\times\mathbb{R}^{d}$ be a
block of axes-parallel boxes which satisfy (\ref{etaadm_cond}). Let $\kappa
\in\left[  0,\kappa_{1}\right[  $ for $\kappa_{1}:=1/6$ and the semi-axes sums
be defined by (\ref{defboldrho}).

\begin{enumerate}
\item[a.] Then, the holomorphic extension of the norm $r$ to $\overrightarrow
{\mathcal{E}}_{\mathbf{a},\mathbf{b},\mathbf{c},\mathbf{d}}\left(
\kappa\right)  $ can be estimated by
\begin{equation}
\left.
\begin{array}
[c]{rl}%
\left\vert r\left(  \mathbf{x},\mathbf{y}\right)  \right\vert \leq & \left(
1+\kappa+2\eta\right)  \operatorname{dist}\left(  \left[  \mathbf{a}%
,\mathbf{b}\right]  ,\left[  \mathbf{c},\mathbf{d}\right]  \right) \\
& \\
\operatorname{Re}r\left(  \mathbf{x},\mathbf{y}\right)  \geq & \left(
1-4\kappa\right)  \operatorname{dist}\left(  \left[  \mathbf{a},\mathbf{b}%
\right]  ,\left[  \mathbf{c},\mathbf{d}\right]  \right)
\end{array}
\right\}  \quad\forall\left(  \mathbf{x},\mathbf{y}\right)  \in\overrightarrow
{\mathcal{E}}_{\mathbf{a},\mathbf{b},\mathbf{c},\mathbf{d}}\left(
\kappa\right)  . \label{upperest}%
\end{equation}

\item[b.] Assume in addition that $\left[  \mathbf{a},\mathbf{b}\right]
\times\left[  \mathbf{c},\mathbf{d}\right]  \subset H_{+}\times H_{-}$. The
function $r_{+}$ (cf. (\ref{Defrrplus})) can be estimated for any $\kappa
\in\left[  0,\min\left\{  \kappa_{1},\frac{1}{3\beta}\right\}  \right[  $ by%
\begin{equation}
\left.
\begin{array}
[c]{rl}%
\left\vert r_{+}\left(  \mathbf{x},\mathbf{y}\right)  \right\vert \leq &
\left(  1+4\beta\right)  \left\vert r\left(  \mathbf{x},\mathbf{y}\right)
\right\vert \\
& \\
\operatorname{Re}r_{+}\left(  \mathbf{x},\mathbf{y}\right)  \geq & \left(
1-3\beta\kappa\right)  \operatorname{Re}r\left(  \mathbf{x},\mathbf{y}\right)
\end{array}
\right\}  \quad\forall\left(  \mathbf{x},\mathbf{y}\right)  \in\overrightarrow
{\mathcal{E}}_{\mathbf{a},\mathbf{b},\mathbf{c},\mathbf{d}}\left(
\kappa\right)  . \label{bothboundsnplus}%
\end{equation}

\end{enumerate}
\end{lemma}%

\proof
\textbf{Part 1: }Proof of (\ref{upperest}).

First, we derive estimates for $\left\Vert \mathbf{x}-\mathbf{y}\right\Vert $.
Again, we only consider the case that the last interval $\left[  a_{d}%
,b_{d}\right]  $ in $\left[  \mathbf{a},\mathbf{b}\right]  $ is extended to
$\mathcal{E}_{a_{d},b_{d}}\left(  \kappa\right)  $. Let $\left(
\mathbf{x},\mathbf{y}\right)  \in\overrightarrow{\mathcal{E}}_{\mathbf{a}%
,\mathbf{b}}^{d}\left(  \kappa\right)  \times\left[  \mathbf{c},\mathbf{d}%
\right]  $ and $\mathfrak{D}$, $\delta$ as in (\ref{deltad}). Let $\xi_{d}%
\in\left[  a_{d},b_{d}\right]  $ and $\mathbf{x}_{\operatorname*{R}}\in\left[
\mathbf{a},\mathbf{b}\right]  $ be as in the previous proof (\ref{defxid}).
The lower estimate follows by combining (\ref{xmylowb}) with (\ref{rests1}):%
\[
\operatorname{Re}r\left(  \mathbf{x},\mathbf{y}\right)  \geq\left(
1-12\kappa^{2}\right)  \left(  1-\kappa\right)  \delta\geq\left(
1-4\kappa\right)  \delta\quad\forall\kappa\in\left[  0,\kappa_{0}\right[  .
\]
The upper estimate is implied by (\ref{xidxd}) and%
\begin{align}
\left\Vert \mathbf{x}-\mathbf{y}\right\Vert  &  \leq\left\Vert \mathbf{x}%
_{\operatorname*{R}}-\mathbf{y}\right\Vert +\left\Vert \mathbf{x}%
_{\operatorname*{R}}-\mathbf{x}\right\Vert \leq\delta+2\mathfrak{D}+\left\vert
\xi_{d}-x_{d}\right\vert \label{normfromabove}\\
&  \leq\delta+2\mathfrak{D}+\kappa\delta\leq\left(  1+2\eta+\kappa\right)
\delta.\nonumber
\end{align}

\textbf{Part 2: }Proof of (\ref{bothboundsnplus}).

To estimate the function $r_{+}$ we consider first the analytic extension for
each of the first $\left(  d-1\right)  $ coordinates. It is straightforward to
obtain from $x_{d}-y_{d}>0$ and $\beta>0:$%
\begin{equation}
\operatorname{Re}r_{+}\left(  \mathbf{x},\mathbf{y}\right)  \geq
\operatorname{Re}r\left(  \mathbf{x},\mathbf{y}\right)  . \label{npluslow1}%
\end{equation}
An upper estimate follows from%
\begin{align}
\left\vert r_{+}\left(  \mathbf{x},\mathbf{y}\right)  \right\vert  &
\leq\left\vert r\left(  \mathbf{x},\mathbf{y}\right)  \right\vert
+\beta\left\Vert \mathbf{x}-\mathbf{y}\right\Vert \label{rplusupperbound}\\
&  \overset{\text{(\ref{rests1})}}{\leq}\left\vert r\left(  \mathbf{x}%
,\mathbf{y}\right)  \right\vert +\frac{\beta}{1-12\kappa^{2}}\operatorname{Re}%
r\left(  \mathbf{x},\mathbf{y}\right)  \leq\left(  1+4\beta\right)  \left\vert
r\left(  \mathbf{x},\mathbf{y}\right)  \right\vert .\nonumber
\end{align}
For $j=d$, the splitting%
\[
x_{d}-y_{d}=\xi_{d}-y_{d}+x_{d}-\xi_{d}%
\]
is employed. By using $\beta>0$, $\xi_{d}-y_{d}>0$, (\ref{xmylowb}), and
(\ref{xidxd}) we obtain%
\begin{align}
\operatorname{Re}r_{+}\left(  \mathbf{x},\mathbf{y}\right)   &  \geq
\operatorname{Re}r\left(  \mathbf{x,y}\right)  +\beta\left(  \xi_{d}%
-y_{d}\right)  -\beta\left\vert x_{d}-\xi_{d}\right\vert \label{npluslow3}\\
&  \geq\operatorname{Re}r\left(  \mathbf{x,y}\right)  -\beta\left\vert
x_{d}-\xi_{d}\right\vert \nonumber\\
&  \geq\operatorname{Re}r\left(  \mathbf{x,y}\right)  -\beta\kappa
\delta.\nonumber
\end{align}
The combination with the second estimate in (\ref{upperest}) yields%
\begin{equation}
\operatorname{Re}r_{+}\left(  \mathbf{x},\mathbf{y}\right)  \geq\left(
1-3\beta\kappa\right)  \operatorname{Re}r\left(  \mathbf{x},\mathbf{y}\right)
. \label{npluslow4}%
\end{equation}
Hence, the lower bound in (\ref{bothboundsnplus}) follows from
(\ref{npluslow1}) and (\ref{npluslow4}). The upper bound follows as for the
extension of the first $\left(  d-1\right)  $ coordinates.%
\endproof

\subsection{Analysis of $\Theta_{\nu,s}^{\operatorname*{illu}}$ and
$\Theta_{\nu,s}^{\operatorname*{refl}}$\label{SecAnaGreenI}}

The proof that $\Theta_{\nu,s}^{\operatorname*{illu}}$ is $\kappa-$slowly
varying follows from the estimates of the Macdonald function and the
holomorphic extension of the norm. The assertion for $\Theta_{\nu
,s}^{\operatorname*{refl}}$ in addition employs a reflection argument.

\begin{lemma}
\label{TheoIlluRefllem}Let $0<\eta\leq\eta_{0}$ be as in Definition
\ref{etaadm}. The families of functions $\mathcal{F}_{\nu}%
^{\operatorname*{illu}}$ and $\mathcal{F}_{\nu}^{\operatorname*{refl}}$ as in
(\ref{slowfamilies}) are $\kappa-$slowly varying: for any $0<\kappa<1/6$, any
$\eta$-admissible block $B=\left[  \mathbf{a},\mathbf{b}\right]  \times\left[
\mathbf{c},\mathbf{d}\right]  \subset H_{+}\times H_{+}$ and $\delta
:=\operatorname{dist}\left(  \left[  \mathbf{a},\mathbf{b}\right]  ,\left[
\mathbf{c},\mathbf{d}\right]  \right)  $ it holds:%
\begin{subequations}
\label{SlowVarIlluRefl}
\end{subequations}%
\begin{align}
\max_{\left(  \mathbf{x},\mathbf{y}\right)  \in\overrightarrow{\mathcal{E}%
}_{\mathbf{a},\mathbf{b},\mathbf{c},\mathbf{d}}\left(  \kappa\right)
}\left\vert \Theta_{\nu,s}^{\operatorname*{illu}}\left(  \mathbf{x}%
,\mathbf{y}\right)  \right\vert  &  \leq C_{\nu}\left(  \frac{\left\vert
s\right\vert }{\delta}\right)  ^{\nu+1/2}M_{\nu+1/2}\left(  \left\vert
s\right\vert \delta\right)  ,\tag{%
\ref{SlowVarIlluRefl}%
a}\label{SlowVarIlluRefla}\\
\max_{\left(  \mathbf{x},\mathbf{y}\right)  \in\overrightarrow{\mathcal{E}%
}_{\mathbf{a},\mathbf{b},\mathbf{c},\mathbf{d}}\left(  \kappa\right)
}\left\vert \Theta_{\nu,s}^{\operatorname*{refl}}\left(  \mathbf{x}%
,\mathbf{y}\right)  \right\vert  &  \leq C_{\nu,\beta}\left(  \frac{\left\vert
s\right\vert }{\delta}\right)  ^{\nu+1/2}M_{\nu+1/2}\left(  \left\vert
s\right\vert \delta\right)  \tag{%
\ref{SlowVarIlluRefl}%
b}\label{SlowVarIlluReflb}%
\end{align}
with $M_{\nu+1/2}$ as in (\ref{majorantdef}).

The constants $\mathfrak{C}_{s}$ and reference functions $\lambda$ in
Definition \ref{DefSlowVar} can be chosen as defined in Theorem
\ref{TheoSlowVar}. In particular, the constants are independent of $\kappa$.
\end{lemma}

%

\proof
First, the claim is proved for $\Theta_{\nu,s}^{\operatorname*{illu}}$
(considered as a function on $\mathbb{R}^{d}\times\mathbb{R}^{d}$). Let
$B:=\left[  \mathbf{a},\mathbf{b}\right]  \times\left[  \mathbf{c}%
,\mathbf{d}\right]  \subset\mathbb{R}^{d}\times\mathbb{R}^{d}$ be an $\eta
$-admissible block. We identify $\left.  \Theta_{\nu,s}^{\operatorname*{illu}%
}\right\vert _{B}$ with its holomorphic extension to $\overrightarrow
{\mathcal{E}}_{\mathbf{a},\mathbf{b},\mathbf{c},\mathbf{d}}\left(
\kappa\right)  $ and consider exemplarily the holomorphic extension to
$\overrightarrow{\mathcal{E}}_{\mathbf{a},\mathbf{b}}^{d}\left(
\kappa\right)  \times\left[  \mathbf{c},\mathbf{d}\right]  $. For $\left(
\mathbf{x},\mathbf{y}\right)  \in\overrightarrow{\mathcal{E}}_{\mathbf{a}%
,\mathbf{b}}^{d}\left(  \kappa\right)  \times\left[  \mathbf{c},\mathbf{d}%
\right]  $ we write $\mathbf{z}:=\mathbf{x}-\mathbf{y}$ and denote by
$r:=r\left(  \mathbf{x},\mathbf{y}\right)  $ the holomorphic norm extension as
in (\ref{defr}). We use the shorthands $\mathfrak{D},\delta$ as in
(\ref{deltad}). Then%
\[
\max_{\left(  \mathbf{x},\mathbf{y}\right)  \in\overrightarrow{\mathcal{E}%
}_{\mathbf{a},\mathbf{b}}^{d}\left(  \kappa\right)  \times\left[
\mathbf{c},\mathbf{d}\right]  }\left\vert \Theta_{\nu,s}^{\operatorname*{illu}%
}\left(  \mathbf{x},\mathbf{y}\right)  \right\vert \leq\frac{1}{\left(
2\pi\right)  ^{\nu+3/2}}\left(  \left\vert \frac{s}{r}\right\vert \right)
^{\nu+1/2}\max_{\left(  \mathbf{x},\mathbf{y}\right)  \in\overrightarrow
{\mathcal{E}}_{\mathbf{a},\mathbf{b}}^{d}\left(  \kappa\right)  \times\left[
\mathbf{c},\mathbf{d}\right]  }\left\vert \operatorname*{e}\nolimits^{sr}%
K_{\nu+1/2}\left(  sr\right)  \right\vert .
\]
From Lemma \ref{Lemrrplus} we conclude that $\operatorname{Re}r>0$ for
$\kappa\in\left[  0,\frac{1}{6}\right[  $ and for $s\in\mathbb{C}_{\geq0}$ we
have%
\[
sr\in\left\{  z\in\mathbb{C}\backslash\mathbb{R}_{\leq0}\mid\frac{\left\vert
s\right\vert \delta}{3}\leq\left\vert z\right\vert \leq\left(  \frac{7}%
{6}+2\eta_{0}\right)  \left\vert s\right\vert \delta\right\}  .
\]
In this way, Corollary \ref{CorFinalMmue} implies%
\[
\max_{\left(  \mathbf{x},\mathbf{y}\right)  \in\overrightarrow{\mathcal{E}%
}_{\mathbf{a},\mathbf{b}}^{d}\left(  \kappa\right)  \times\left[
\mathbf{c},\mathbf{d}\right]  }\left\vert \operatorname*{e}\nolimits^{sr}%
K_{\nu+1/2}\left(  sr\right)  \right\vert \leq C_{\nu}M_{\nu+1/2}\left(
\left\vert sr\right\vert \right)
\]
and%
\[
\max_{\left(  \mathbf{x},\mathbf{y}\right)  \in\overrightarrow{\mathcal{E}%
}_{\mathbf{a},\mathbf{b}}^{d}\left(  \kappa\right)  \times\left[
\mathbf{c},\mathbf{d}\right]  }\left\vert \Theta_{\nu,s}^{\operatorname*{illu}%
}\left(  \mathbf{x},\mathbf{y}\right)  \right\vert \leq C_{\nu}\left(
\frac{\left\vert s\right\vert }{\delta}\right)  ^{\nu+1/2}M_{\nu+1/2}\left(
\left\vert s\right\vert \delta\right)  ,
\]
where $C_{\nu}$ only depends on $\nu$.\medskip

To show that $\Theta_{\nu,s}^{\operatorname*{refl}}$ is $\kappa-$slowly
varying we define the reflected box%
\[
\mathbf{R}\left[  \mathbf{c},\mathbf{d}\right]  :=\left\{  \mathbf{Ry}%
:\mathbf{y}\in\left[  \mathbf{c,d}\right]  \right\}  \subset H_{-}%
\]
and set%
\[
\mathbf{z}_{-}:=\mathbf{x}-\mathbf{Ry,}\text{\quad}d_{-}:=\sqrt{\left\langle
\mathbf{z}_{-},\mathbf{z}_{-}\right\rangle },\quad z_{d}^{-}:=x_{d}+y_{d}.
\]
The $\eta$-admissibility of $\left[  \mathbf{a},\mathbf{b}\right]
\times\left[  \mathbf{c},\mathbf{d}\right]  $ implies the $\eta$-admissibility
of $\left[  \mathbf{a},\mathbf{b}\right]  \times\mathbf{R}\left[
\mathbf{c},\mathbf{d}\right]  $ as can be seen from the following reasoning:
choose $\left(  \mathbf{x}_{0},\mathbf{y}_{0}\right)  \in\left[
\mathbf{a},\mathbf{b}\right]  \times\left[  \mathbf{c},\mathbf{d}\right]  $
such that
\[
\delta_{\operatorname*{refl}}:=\operatorname{dist}\left(  \left[
\mathbf{a},\mathbf{b}\right]  ,\mathbf{R}\left[  \mathbf{c},\mathbf{d}\right]
\right)  =\left\Vert \mathbf{x}_{0}-\mathbf{Ry}_{0}\right\Vert .
\]
Then, $\delta=\operatorname{dist}\left(  \left[  \mathbf{a},\mathbf{b}\right]
,\left[  \mathbf{c},\mathbf{d}\right]  \right)  $ can be estimated by using
$x_{0,d}-y_{0,d}\leq x_{0,d}+y_{0,d}$ for any $\mathbf{x}_{0},\mathbf{y}%
_{0}\in H_{+}:$%
\begin{equation}
\delta\leq\left\Vert \mathbf{x}_{0}-\mathbf{y}_{0}\right\Vert =\sqrt
{\left\Vert \mathbf{x}_{0}^{\prime}-\mathbf{y}_{0}^{\prime}\right\Vert
^{2}+\left(  x_{0,d}-y_{0,d}\right)  ^{2}}\leq\left\Vert \mathbf{x}%
_{0}-\mathbf{Ry}_{0}\right\Vert =\delta_{\operatorname*{refl}} \label{drefld}%
\end{equation}
and the auxiliary statement follows from%
\begin{equation}
\max\left\{  \operatorname*{diam}\left[  \mathbf{a},\mathbf{b}\right]
,\operatorname*{diam}\mathbf{R}\left[  \mathbf{c},\mathbf{d}\right]  \right\}
=\max\left\{  \operatorname*{diam}\left[  \mathbf{a},\mathbf{b}\right]
,\operatorname*{diam}\left[  \mathbf{c},\mathbf{d}\right]  \right\}  \leq
\eta\delta\leq\eta\delta_{\operatorname*{refl}}. \label{admrefl}%
\end{equation}

Recall that $\Theta_{\nu,s}^{\operatorname*{refl}}$ is defined in
(\ref{slowGreenb}) via the function $\sigma_{\nu}$ from (\ref{Defpsinues}) and
we estimate the prefactor in (\ref{Defpsinues}) first. By the same arguments
as for (\ref{npluslow4}) it follows that%
\begin{align*}
\operatorname{Re}\left(  \beta d_{-}+z_{d}^{-}\right)   &  \geq\beta\left(
1-3\kappa\right)  \operatorname{Re}d_{-}\geq\frac{\beta}{2}\operatorname{Re}%
d_{-}\\
&  \overset{\text{(\ref{rests1})}}{\geq}\frac{\beta}{2}\left(  1-12\kappa
^{2}\right)  \left\Vert \mathbf{z}_{-}\right\Vert \overset{\kappa<1/6}{\geq
}\frac{\beta}{3}\left\Vert \mathbf{z}_{-}\right\Vert \qquad\forall\left(
\mathbf{x},\mathbf{y}\right)  \in\overrightarrow{\mathcal{E}}_{\mathbf{a}%
,\mathbf{b},\mathbf{c},\mathbf{d}}\left(  \kappa\right)  .
\end{align*}
For the numerator in the prefactor of $\sigma_{\nu}$ it holds%
\[
\left\vert z_{d}^{-}-\beta d_{-}\right\vert \leq\left(  1+\beta\right)
\left\Vert \mathbf{z}_{-}\right\Vert
\]
so that%
\[
\max_{\left(  \mathbf{x},\mathbf{y}\right)  \in\overrightarrow{\mathcal{E}%
}_{\mathbf{a},\mathbf{b},\mathbf{c},\mathbf{d}}\left(  \kappa\right)
}\left\vert \frac{z_{d}^{-}-\beta d_{-}}{z_{d}^{-}+\beta d_{-}}\right\vert
\leq3\frac{1+\beta}{\beta}.
\]
The estimates for the function $\Theta_{\nu,s}^{\operatorname*{illu}}$ then
imply%
\begin{align*}
\max_{\left(  \mathbf{x},\mathbf{y}\right)  \in\overrightarrow{\mathcal{E}%
}_{\mathbf{a},\mathbf{b}}^{d}\left(  \kappa\right)  \times\left[
\mathbf{c},\mathbf{d}\right]  }\left\vert \Theta_{\nu,s}^{\operatorname*{refl}%
}\left(  \mathbf{x},\mathbf{y}\right)  \right\vert  &  =3\frac{1+\beta}{\beta
}\max_{\left(  \mathbf{x},\mathbf{y}\right)  \in\overrightarrow{\mathcal{E}%
}_{\mathbf{a},\mathbf{b}}^{d}\left(  \kappa\right)  \times\mathbf{R}\left[
\mathbf{c},\mathbf{d}\right]  }\left\vert \Theta_{\nu,s}^{\operatorname*{illu}%
}\left(  \mathbf{x},\mathbf{y}\right)  \right\vert \\
&  \leq C_{\nu,\beta}\left(  \frac{\left\vert s\right\vert }{\delta
_{\operatorname*{refl}}}\right)  ^{\nu+1/2}M_{\nu+1/2}\left(  \left\vert
s\right\vert \delta_{\operatorname*{refl}}\right)  .
\end{align*}
Since $\delta\leq\delta_{\operatorname*{refl}}$ and $M_{\nu+1/2}$ is strictly
decreasing, the assertion follows by adjusting $C_{\nu,\beta}$.
\endproof

\subsection{Analysis of $\Theta_{\nu,s}^{\operatorname*{imp}}$%
\label{SecAnaGreenII}}

In this section it will be shown that the family of functions $\mathcal{F}%
_{\nu}^{\operatorname*{imp}}$ (see (\ref{slowfamiliesc})) is $\kappa-$slowly
varying. We employ the definition (\ref{slowGreenc}) and the definition of
$\psi_{\nu,s}$ as in (\ref{Defpsinues})%
\begin{equation}
\Theta_{\nu,s}^{\operatorname*{imp}}\left(  \mathbf{x},\mathbf{y}\right)
=-\frac{\beta}{\pi}\left(  \frac{s^{2}}{2\pi}\right)  ^{\nu+1/2}\frac{1}%
{s}\int_{0}^{\infty}\operatorname*{e}\nolimits^{-sy}q_{\nu}\left(
\mathbf{z},y\right)  dy. \label{repsTheta}%
\end{equation}
and%
\[
\mathbf{z}:=\mathbf{x}-\mathbf{Ry\quad}\text{and\quad}r:=\left\Vert
\mathbf{z}\right\Vert .
\]

We start with an estimate of the function $q_{\nu}\left(  \mathbf{z},y\right)
$ as in (\ref{defqnue}). We employ the notation and setting: let $B=\left[
\mathbf{a},\mathbf{b}\right]  \times\left[  \mathbf{c},\mathbf{d}\right]
\subset H_{+}\times H_{+}$ be an $\eta$-admissible block and $\delta
:=\operatorname{dist}\left(  \left[  \mathbf{a},\mathbf{b}\right]  ,\left[
\mathbf{c},\mathbf{d}\right]  \right)  $. For $\left(  \mathbf{x}%
,\mathbf{y}\right)  \in\overrightarrow{\mathcal{E}}_{\mathbf{a},\mathbf{b}%
,\mathbf{c},\mathbf{d}}\left(  \kappa\right)  $ let $\mathbf{z}:=\mathbf{x}%
-\mathbf{Ry}$. The holomorphic functions $r,r_{+}:\overrightarrow{\mathcal{E}%
}_{\mathbf{a},\mathbf{b},\mathbf{c},\mathbf{d}}\left(  \kappa\right)
\rightarrow\mathbb{C}$ are given by%
\[
r\left(  \mathbf{z}\right)  =\sqrt{\left\langle \mathbf{z},\mathbf{z}%
\right\rangle }\quad\text{and}\quad r_{+}\left(  \mathbf{z}\right)  =r\left(
\mathbf{z}\right)  +\beta z_{d}\text{.}%
\]
Recall the definition of the function $\tilde{\mu}$ as in (\ref{defsmue}). An
explicit calculation yields $q_{\nu}=q_{\nu}^{\operatorname*{I}}+q_{\nu
}^{\operatorname*{II}}$ with%
\begin{align*}
q_{\nu}^{\operatorname*{I}}\left(  \mathbf{z},y\right)   &  :=s\left(
\tilde{\mu}^{\prime}-\frac{\tilde{\mu}\left(  \tilde{\mu}+\beta t\right)
}{\left(  t+\beta\tilde{\mu}\right)  ^{2}}\right)  \frac{\operatorname*{e}%
\nolimits^{s\tilde{\mu}}K_{\nu+1/2}\left(  s\tilde{\mu}\right)  }{\left(
t+\beta\tilde{\mu}\right)  \left(  s\tilde{\mu}\right)  ^{\nu+1/2}},\\
q_{\nu}^{\operatorname*{II}}\left(  \mathbf{z},y\right)   &  :=s^{2}\tilde
{\mu}^{\prime}\frac{\tilde{\mu}}{t+\beta\tilde{\mu}}\left(  \frac
{\operatorname*{e}\nolimits^{s\tilde{\mu}}\left(  K_{\nu+1/2}\left(
s\tilde{\mu}\right)  -K_{\nu+3/2}\left(  s\tilde{\mu}\right)  \right)
}{\left(  s\tilde{\mu}\right)  ^{\nu+1/2}}\right)  .
\end{align*}
Clearly, a key role in the analysis of $q_{\nu}$ and, in turn, of $\Theta
_{\nu,s}^{\operatorname*{imp}}$ is played by the function $\tilde{\mu}$ and
two-sided estimated are stated in the following lemma.

\begin{lemma}
\label{Lemmueestgen}Let $\left(  \left[  \mathbf{a},\mathbf{b}\right]
\times\left[  \mathbf{c},\mathbf{d}\right]  \right)  \subset H_{+}\times
H_{+}$ be a block of axes-parallel boxes which satisfy (\ref{etaadm_cond}).
There exist numbers $C_{\mathcal{E}}\geq6$, $C_{0}>1$, and $c_{2}>0$
independent of all parameters and functions such that for any $\kappa
\in\left[  0,\frac{\beta^{2}}{C_{\mathcal{E}}\left(  1+\beta\right)  ^{3}%
}\right[  $ and semi-axes sums defined by (\ref{defboldrho}) and any $\left(
\mathbf{x},\mathbf{y}\right)  \in\overrightarrow{\mathcal{E}}_{\mathbf{a}%
,\mathbf{b},\mathbf{c},\mathbf{d}}\left(  \kappa\right)  $ it holds for
$\mathbf{z}:=\mathbf{x}-\mathbf{Ry}$ and $y\geq0:$%
\begin{subequations}
\label{mueest}
\end{subequations}%
\begin{equation}
\left\vert \tilde{\mu}\left(  \mathbf{z},y\right)  \right\vert \leq C_{0}%
\frac{1+\beta}{\beta}\left(  y+\left(  1+\beta\right)  \left\Vert
\mathbf{z}\right\Vert \right)  \tag{%
\ref{mueest}%
a}\label{mueesta}%
\end{equation}
and%
\begin{equation}
\operatorname{Re}\tilde{\mu}\left(  \mathbf{z},y\right)  \geq C_{0}^{-1}%
\frac{\beta^{2}}{\left(  1+\beta\right)  ^{4}}\left(  y+\left\Vert
\mathbf{z}\right\Vert \right)  . \tag{%
\ref{mueest}%
b}\label{mueestb}%
\end{equation}
For the sum $t+\beta\tilde{\mu}$ it holds%
\begin{equation}
\operatorname{Re}\left(  t\left(  \mathbf{z},y\right)  +\beta\tilde{\mu
}\left(  \mathbf{z},y\right)  \right)  \geq c_{2}\left(  y+\beta\left\Vert
\mathbf{z}\right\Vert \right)  \tag{%
\ref{mueest}%
c}\label{mueestb2}%
\end{equation}
and for the derivative:%
\begin{equation}
\left\vert \frac{\partial\tilde{\mu}\left(  \mathbf{z},y\right)  }{\partial
y}\right\vert \leq C_{0}\left(  \frac{1+\beta}{\beta}\right)  ^{2}. \tag{%
\ref{mueest}%
d}\label{mueestc}%
\end{equation}

\end{lemma}

\begin{remark}
The upper bound for the range of $\kappa$ in Lemma \ref{Lemmueestgen}, i.e.,
$\frac{\beta^{2}}{C_{\mathcal{E}}\left(  1+\beta\right)  ^{3}}$, takes its
maximal value $\frac{4}{27}C_{\mathcal{E}}^{-1}$ at $\beta=2$; hence, the
upper bound for $\kappa$ tends linearly to zero as $C_{\mathcal{E}}$
increases. It is also straightforward to verify that $C_{\mathcal{E}}>6$
implies
\[
\frac{\beta^{2}}{C_{\mathcal{E}}\left(  1+\beta\right)  ^{3}}\leq\frac{1}%
{6}\min\left\{  1,\beta^{-1}\right\}  .
\]

\end{remark}

The proof of Lemma \ref{Lemmueestgen} is fairly technical and postponed to
Section \ref{Secmue}.

Lemma \ref{Lemmueestgen} allows us to estimate the function $q_{\nu}$.

\begin{lemma}
\label{LemVarimpII}Let $0<\eta\leq\eta_{0}$ be as in Definition \ref{etaadm}.
Let $\left(  \left[  \mathbf{a},\mathbf{b}\right]  ,\left[  \mathbf{c}%
,\mathbf{d}\right]  \right)  \subset H_{+}\times H_{+}$ be a block of
axes-parallel boxes which satisfy (\ref{etaadm_cond}). Let $C_{\mathcal{E}%
},C_{0}>0$ be as in Lemma \ref{Lemmueestgen}.

\begin{enumerate}
\item For any $\kappa\in\left[  0,\frac{\beta^{2}}{C_{\mathcal{E}}\left(
1+\beta\right)  ^{3}}\right[  $ and semi-axes sums defined by
(\ref{defboldrho}), any $\left(  \mathbf{x},\mathbf{y}\right)  \in
\overrightarrow{\mathcal{E}}_{\mathbf{a},\mathbf{b},\mathbf{c},\mathbf{d}%
}\left(  \kappa\right)  $ and $\mathbf{z}:=\mathbf{x}-\mathbf{Ry}$, the
integral $\int_{0}^{\infty}q_{\nu}\left(  \mathbf{z},y\right)  dy$ exists as
an improper Riemann integral.

\item The family of functions $\mathcal{F}_{\nu}^{\operatorname*{imp}}$ in
(\ref{slowfamiliesc}) is $\kappa-$slowly varying: for any $\kappa\in\left[
0,\frac{\beta^{2}}{C_{\mathcal{E}}\left(  1+\beta\right)  ^{3}}\right[  $, any
$\eta$-admissible block $B=\left[  \mathbf{a},\mathbf{b}\right]  \times\left[
\mathbf{c},\mathbf{d}\right]  \subset H_{+}\times H_{+}$ and $\delta
:=\operatorname{dist}\left(  \left[  \mathbf{a},\mathbf{b}\right]  ,\left[
\mathbf{c},\mathbf{d}\right]  \right)  $ it holds
\[
\max_{\left(  \mathbf{x},\mathbf{y}\right)  \in\overrightarrow{\mathcal{E}%
}_{\mathbf{a},\mathbf{b},\mathbf{c},\mathbf{d}}\left(
\mbox{\boldmath$ \rho$}%
_{1},%
\mbox{\boldmath$ \rho$}%
_{2}\right)  }\left\vert \Theta_{\nu,s}^{\operatorname*{imp}}\left(
\mathbf{x},\mathbf{y}\right)  \right\vert \leq C_{\nu,\beta}\left(
\frac{\left\vert s\right\vert }{\delta}\right)  ^{\nu+1/2}W_{\nu+1/2}\left(
\left\vert s\right\vert \delta\right)
\]
with $W_{\mu}$ as in (\ref{DefMWb}). The constant $\mathfrak{C}_{s}$ and
reference function $\lambda$ in Definition \ref{DefSlowVar} can be chosen as
defined in Theorem \ref{TheoSlowVar}. In particular, the constants are
independent of $\kappa$.
\end{enumerate}
\end{lemma}

%

\proof
Let $\left(  \mathbf{x},\mathbf{y}\right)  \in\overrightarrow{\mathcal{E}%
}_{\mathbf{a},\mathbf{b},\mathbf{c},\mathbf{d}}\left(  \kappa\right)  $ and
$\mathbf{z}:=\mathbf{x}-\mathbf{Ry}$. We use the shorthands
$r_{\operatorname*{R}}=\left\Vert \mathbf{z}\right\Vert $ and $\tilde{\mu
}=\tilde{\mu}\left(  \mathbf{z},y\right)  $. The combination of
(\ref{bothboundsnplus}) with (\ref{rests1}) leads to%
\begin{equation}
\operatorname{Re}r_{+}\left(  \mathbf{z}\right)  \geq\left(  1-3\beta
\kappa\right)  \left(  1-12\kappa^{2}\right)  r_{\operatorname*{R}}\geq
\frac{1}{3}r_{\operatorname*{R}}\quad\text{and\quad}\operatorname{Re}\left(
r_{+}+y\right)  \geq y+\frac{r_{\operatorname*{R}}}{3}.
\label{estypluerplusbelow}%
\end{equation}
Then the inequalities (\ref{mueest}) and the relation (cf. (\ref{defyvont}))
$\tilde{\mu}+\beta t=y+r_{+}$ imply
\begin{align*}
\left\vert \tilde{\mu}^{\prime}-\frac{\tilde{\mu}\left(  \tilde{\mu}+\beta
t\right)  }{\left(  t+\beta\tilde{\mu}\right)  ^{2}}\right\vert  &  \leq
C_{0}\left(  \frac{1+\beta}{\beta}\right)  ^{2}+\frac{C_{0}}{c_{2}^{2}}%
\frac{1+\beta}{\beta}\left(  \frac{y+\left(  1+\beta\right)
r_{\operatorname*{R}}}{y+\beta r_{\operatorname*{R}}}\right)  ^{2}\\
&  \leq C_{0}\left(  \frac{1+\beta}{\beta}\right)  ^{2}+\frac{C_{0}}{c_{2}%
^{2}}\frac{\left(  1+\beta\right)  ^{3}}{\beta^{3}}\leq C_{0}^{\prime}%
\frac{\left(  1+\beta\right)  ^{3}}{\beta^{3}}.
\end{align*}
From Lemma \ref{Lemmueestgen} and $s\in\mathbb{C}_{\geq0}$ it follows%
\[
s\tilde{\mu}\in\left\{  \zeta\in\mathbb{C}\backslash\mathbb{R}_{\leq0}\mid
C_{0}^{-1}\frac{\beta^{2}}{\left(  1+\beta\right)  ^{4}}\left\vert
s\right\vert \left(  y+r_{\operatorname*{R}}\right)  \leq\left\vert
\zeta\right\vert \leq C_{0}\frac{1+\beta}{\beta}\left\vert s\right\vert
\left(  y+\left(  1+\beta\right)  r_{\operatorname*{R}}\right)  \right\}  .
\]
This allows us to apply Corollary \ref{CorFinalMmue} so that%
\begin{align*}
\left\vert q_{\nu}^{\operatorname*{I}}\left(  \mathbf{z},y\right)
\right\vert  &  \leq C_{0}^{\prime}\frac{\left(  1+\beta\right)  ^{3}}%
{\beta^{3}}\left\vert s\right\vert \frac{M_{\nu+1/2}\left(  C_{0}^{-1}%
\frac{\beta^{2}}{\left(  1+\beta\right)  ^{4}}\left\vert s\right\vert \left(
y+r_{\operatorname*{R}}\right)  \right)  }{c_{2}\left(  y+\beta
r_{\operatorname*{R}}\right)  \left(  C_{0}^{-1}\frac{\beta^{2}}{\left(
1+\beta\right)  ^{4}}\left\vert s\right\vert \left(  y+r_{\operatorname*{R}%
}\right)  \right)  ^{\nu+1/2}}\\
&  \leq C_{\nu,\beta}\left\vert s\right\vert ^{2}\frac{M_{\nu+1/2}\left(
C_{0}^{-1}\frac{\beta^{2}}{\left(  1+\beta\right)  ^{4}}\left\vert
s\right\vert \left(  y+r_{\operatorname*{R}}\right)  \right)  }{\left(
\left\vert s\right\vert \left(  y+r_{\operatorname*{R}}\right)  \right)
^{\nu+3/2}}.
\end{align*}
To estimate $q_{\nu}^{\operatorname*{II}}$, we combine Lemma
\ref{Lemmueestgen}, Corollary \ref{CorFinalMmue}, and
(\ref{estypluerplusbelow}) and obtain%
\[
\left\vert q_{\nu}^{\operatorname*{II}}\left(  \mathbf{z},y\right)
\right\vert \leq C_{\nu,\beta}\left\vert s\right\vert ^{2}\frac{N_{\nu
+1/2}\left(  C_{0}^{-1}\frac{\beta^{2}}{\left(  1+\beta\right)  ^{4}%
}\left\vert s\right\vert \left(  y+r_{\operatorname*{R}}\right)  \right)
}{\left(  \left\vert s\right\vert \left(  y+r_{\operatorname*{R}}\right)
\right)  ^{\nu+3/2}}.
\]

It follows that a majorant for the function $q_{\nu}\left(  \mathbf{z}%
,\cdot\right)  $ is given by a rational function without poles for
$y\in\left]  0,\infty\right[  $ since $r_{R}>0$. For large $z$ we use that
$M_{\nu+1/2}\left(  z\right)  $ as well as $N_{\nu+1/2}\left(  z\right)  $
decay with a speed of $\left\vert z\right\vert ^{-1/2}$
(cf.\ Corollary~\ref{CorFinalMmue}) so that $q_{\nu}\left(  \mathbf{z}%
,\cdot\right)  $ decays with a speed of $O\left(  \left(
y+r_{\operatorname*{R}}\right)  ^{-\nu-2}\right)  $. From $\nu\geq-1/2$ we
conclude that $\int_{0}^{\infty}q_{\nu}\left(  \mathbf{z},y\right)  dy$ exists
as an improper Riemann integral.

For the estimate of the function $\Theta_{\nu,s}^{\operatorname*{imp}}$ we
recall that $M_{\nu+1/2}$ is a majorant of $N_{\nu+1/2}$ (see (\ref{NmueMmue}%
)) so that%
\[
\left\vert q_{\nu}\left(  \mathbf{z},y\right)  \right\vert \leq C_{\nu,\beta
}\left\vert s\right\vert ^{2}\frac{M_{\nu+1/2}\left(  C_{0}^{-1}\frac
{\beta^{2}}{\left(  1+\beta\right)  ^{4}}\left\vert s\right\vert \left(
y+r_{\operatorname*{R}}\right)  \right)  }{\left(  \left\vert s\right\vert
\left(  y+r_{\operatorname*{R}}\right)  \right)  ^{\nu+3/2}}.
\]
This leads to%
\begin{equation}
\left\vert \Theta_{\nu,s}^{\operatorname*{imp}}\left(  \mathbf{x}%
,\mathbf{y}\right)  \right\vert \leq C_{\nu,\beta}\left\vert s\right\vert
^{2\left(  \nu+1\right)  }\int_{0}^{\infty}\frac{M_{\nu+1/2}\left(  C_{0}%
^{-1}\frac{\beta^{2}}{\left(  1+\beta\right)  ^{4}}\left\vert s\right\vert
\left(  y+r_{\operatorname*{R}}\right)  \right)  }{\left(  \left\vert
s\right\vert \left(  y+r_{\operatorname*{R}}\right)  \right)  ^{\nu+3/2}}dy.
\label{thetanuesII}%
\end{equation}
In view of the piecewise definition of $M_{\mu}$ in (\ref{majorantdef}) we set%
\begin{equation}
b:=\left\{
\begin{array}
[c]{ll}%
1/\left\vert s\right\vert -r_{\operatorname*{R}} & \text{for }\left\vert
s\right\vert r_{\operatorname*{R}}<1,\\
0 & \text{for }\left\vert s\right\vert r_{\operatorname*{R}}\geq1
\end{array}
\right.  \label{defb}%
\end{equation}
and observe%
\[
\left\vert s\right\vert \left(  b+r_{\operatorname*{R}}\right)  =\max\left\{
1,\left\vert s\right\vert r_{\operatorname*{R}}\right\}  .
\]
This leads to the splitting $\left(  \int_{0}^{b}\ldots+\int_{b}^{\infty
}\ldots\right)  $ of the integral (\ref{thetanuesII}). For the second summand
we get%
\begin{align*}
\int_{b}^{\infty}\frac{M_{\nu+1/2}\left(  C_{0}^{-1}\frac{\beta^{2}}{\left(
1+\beta\right)  ^{4}}\left\vert s\right\vert \left(  y+r_{\operatorname*{R}%
}\right)  \right)  }{\left(  \left\vert s\right\vert \left(
y+r_{\operatorname*{R}}\right)  \right)  ^{\nu+3/2}}dy  &  =C_{\nu,\beta}%
\int_{b}^{\infty}\left(  \left\vert s\right\vert \left(
y+r_{\operatorname*{R}}\right)  \right)  ^{-\left(  \nu+2\right)  }dy\\
&  =C_{\nu,\beta}\left\vert s\right\vert ^{-1}\left(  \left\vert s\right\vert
\left(  b+r_{\operatorname*{R}}\right)  \right)  ^{-\left(  \nu+1\right)  }\\
&  \leq C_{\nu,\beta}\frac{1}{\left\vert s\right\vert \left(  \max\left\{
1,\left\vert s\right\vert r_{\operatorname*{R}}\right\}  \right)  ^{\nu+1}}.
\end{align*}
The integral over $\left(  0,b\right)  $ is non-zero only for the first case
in (\ref{defb}) so that $\left\vert s\right\vert r_{\operatorname*{R}}<1$. We
first consider the case $\nu\in\mathbb{Z}_{\geq0}^{\operatorname*{half}}$.
Then%
\begin{align*}
\int_{0}^{b}\frac{M_{\nu+1/2}\left(  C_{0}^{-1}\frac{\beta^{2}}{\left(
1+\beta\right)  ^{4}}\left\vert s\right\vert \left(  y+r_{\operatorname*{R}%
}\right)  \right)  }{\left(  \left\vert s\right\vert \left(
y+r_{\operatorname*{R}}\right)  \right)  ^{\nu+3/2}}dy  &  =C_{\nu,\beta}%
\int_{0}^{1/\left\vert s\right\vert -r_{\operatorname*{R}}}\frac{1}{\left(
\left\vert s\right\vert \left(  y+r_{\operatorname*{R}}\right)  \right)
^{2\nu+2}}dy\\
&  \leq C_{\nu,\beta}\frac{1}{\left\vert s\right\vert }\frac{1}{\left(
\left\vert s\right\vert r_{\operatorname*{R}}\right)  ^{2\nu+1}}.
\end{align*}
The estimate for the full integral is then given by%
\[
\int_{0}^{\infty}\frac{M_{\nu+1/2}\left(  C_{0}^{-1}\frac{\beta^{2}}{\left(
1+\beta\right)  ^{4}}\left\vert s\right\vert \left(  y+r_{\operatorname*{R}%
}\right)  \right)  }{\left(  \left\vert s\right\vert \left(
y+r_{\operatorname*{R}}\right)  \right)  ^{\nu+3/2}}dy\leq C_{\nu,\beta
}\left\{
\begin{array}
[c]{lc}%
\frac{1}{\left\vert s\right\vert \left(  \left\vert s\right\vert
r_{\operatorname*{R}}\right)  ^{\nu+1}} & \text{for }\left\vert s\right\vert
r_{\operatorname*{R}}\geq1,\\
\frac{1}{\left\vert s\right\vert \left(  \left\vert s\right\vert
r_{\operatorname*{R}}\right)  ^{2\nu+1}} & \text{for }\left\vert s\right\vert
r_{\operatorname*{R}}<1.
\end{array}
\right.
\]
It remains to consider the integral over $\left(  0,b\right)  $ for $\nu=-1/2$
and $\left\vert s\right\vert r_{\operatorname*{R}}<1$. We get\footnote{It
holds%
\[
\int_{0}^{1/\left\vert s\right\vert -r_{\operatorname*{R}}}\frac{\left\vert
\ln\left(  \left\vert s\right\vert \left(  y+r_{\operatorname*{R}}\right)
\right)  \right\vert }{\left\vert s\right\vert \left(  y+r_{\operatorname*{R}%
}\right)  }dy=\frac{\ln^{2}\left(  \left\vert s\right\vert
r_{\operatorname*{R}}\right)  }{2|s|}%
\]
}%
\begin{align*}
\int_{0}^{b}\frac{M_{0}\left(  C_{0}^{-1}\frac{\beta^{2}}{\left(
1+\beta\right)  ^{4}}\left\vert s\right\vert \left(  y+r_{\operatorname*{R}%
}\right)  \right)  }{\left\vert s\right\vert \left(  y+r_{\operatorname*{R}%
}\right)  }dy  &  \leq C_{-1/2,\beta}\int_{0}^{1/\left\vert s\right\vert
-r_{\operatorname*{R}}}\frac{\left(  1+\left\vert \ln\left(  \left\vert
s\right\vert \left(  y+r_{\operatorname*{R}}\right)  \right)  \right\vert
\right)  }{\left\vert s\right\vert \left(  y+r_{\operatorname*{R}}\right)
}dy\\
&  =C_{-1/2,\beta}\frac{\left\vert \ln\left(  \left\vert s\right\vert
r_{\operatorname*{R}}\right)  \right\vert }{|s|}+C_{-1/2,\beta}\left(
\frac{\ln^{2}\left(  \left\vert s\right\vert r_{\operatorname*{R}}\right)
}{2|s|}\right) \\
&  \leq C_{-1/2,\beta}\frac{1+\ln^{2}\left(  \left\vert s\right\vert
r_{\operatorname*{R}}\right)  }{\left\vert s\right\vert }.
\end{align*}
This leads to%
\begin{equation}
\left\vert \Theta_{\nu,s}^{\operatorname*{imp}}\left(  \mathbf{x}%
,\mathbf{y}\right)  \right\vert \leq C_{\nu,\beta}\left\{
\begin{array}
[c]{ll}%
\left(  \frac{\left\vert s\right\vert }{r_{\operatorname*{R}}}\right)  ^{\nu
}\frac{1}{r_{\operatorname*{R}}} & \text{for }\left\vert s\right\vert
r_{\operatorname*{R}}\geq1,\\
\frac{1}{r_{\operatorname*{R}}^{2\nu+1}} & \text{for }\left\vert s\right\vert
r_{\operatorname*{R}}<1\wedge\nu\in\mathbb{Z}_{\geq0}^{\operatorname*{half}%
},\\
1+\ln^{2}\left(  \left\vert s\right\vert r_{\operatorname*{R}}\right)  &
\text{for }\left\vert s\right\vert r_{\operatorname*{R}}<1\wedge\nu=-1/2
\end{array}
\right.  \label{thetaimp}%
\end{equation}
from which the assertion follows.

\begin{remark}
Let $\nu\in\mathbb{Z}_{\geq0}^{\operatorname*{half}}$. Estimate
(\ref{thetaimp}) shows that the dominant singular behaviour of $\Theta_{\nu
,s}^{\operatorname*{imp}}\left(  \mathbf{x},\mathbf{y}\right)  $ can be
estimated for small $\left\Vert \mathbf{x}-\mathbf{Ry}\right\Vert $ by
$O\left(  \left\Vert \mathbf{x}-\mathbf{Ry}\right\Vert ^{-2\nu-1}\right)  $.
We emphasize that a singular expansion of $\left\vert \Theta_{\nu
,s}^{\operatorname*{imp}}\left(  \mathbf{x},\mathbf{y}\right)  \right\vert $
for small $\left\Vert \mathbf{x}-\mathbf{Ry}\right\Vert $ may also contain
lower order terms with logarithmic singularities. Indeed, for $\nu=0$ and
$\beta=1$ the representation\footnote{A similar representation can be found
for different impedance parameter $\beta$ in \cite[(5.105)]{hoernig2010green}%
.} of $G_{\operatorname*{half}}$%
\begin{equation}
G_{\operatorname*{half}}\left(  \mathbf{x},\mathbf{y}\right)  :=\frac
{\operatorname*{e}\nolimits^{-s\left\Vert \mathbf{x}-\mathbf{y}\right\Vert }%
}{4\pi\left\Vert \mathbf{x}-\mathbf{y}\right\Vert }+\frac{\operatorname*{e}%
\nolimits^{-s\left\Vert \mathbf{x}-\mathbf{Ry}\right\Vert }}{4\pi\left\Vert
\mathbf{x}-\mathbf{Ry}\right\Vert }-\frac{s\operatorname*{e}%
\nolimits^{-s\left\Vert \mathbf{z}\right\Vert }}{2\pi}U\left(  1,1,s\left(
\left\Vert \mathbf{z}\right\Vert +z_{3}\right)  \right)  \label{f}%
\end{equation}
with Tricomi's (confluent hypergeometric) function $U\left(  a,b,z\right)  $
(see \cite[13.2.6]{NIST:DLMF}) is derived in \cite[Sec. 4]%
{LinMelenkSauter_Gimp_I}. Clearly, the dominant singular behaviour of the
first two summands for small arguments is $\left\Vert \mathbf{x}%
-\mathbf{y}\right\Vert ^{-1}$ for the first and $\left\Vert \mathbf{x}%
-\mathbf{Ry}\right\Vert ^{-1}$ for the second summand which shows that
(\ref{thetaimp}) is sharp in this case. However, the dominant singular
behaviour of Tricomi's function is logarithmic: $U\left(  1,1,\rho\right)
=-\left(  \gamma+1+\log\rho\right)  +o\left(  1\right)  $ as $\rho
\rightarrow0$ (see \cite[Chap. 48:9]{OldhamAtlas}) and the third term in
(\ref{f}) has a singularity of lower order compared to the second summand.
\end{remark}

\section{Estimate of $\tilde{\mu}$\label{Secmue}}

We introduce the notation%
\[
\chi=a^{2}-\left(  1-\beta^{2}\right)  \omega^{2}\quad\text{with\quad
}a:=y+r_{+},\quad\omega^{2}:=\left\langle \mathbf{z}^{\prime},\mathbf{z}%
^{\prime}\right\rangle
\]
and obtain the compact representation of $\tilde{\mu}$ and $t$:
\begin{equation}
\tilde{\mu}\left(  \mathbf{z},t\right)  =\frac{a\sqrt{\chi}+\beta\omega^{2}%
}{\beta a+\sqrt{\chi}}\quad\text{and\quad}t\left(  \mathbf{z},y\right)
=\frac{a^{2}-\omega^{2}}{\beta a+\sqrt{\chi}}. \label{muecomprep}%
\end{equation}

\textbf{Proof of Lemma \ref{Lemmueestgen}. }First, estimates for $\chi$ will
be derived. An explicit computation yields%
\begin{align}
\chi &  =\left(  y+r+\beta z_{d}\right)  ^{2}-\left(  1-\beta^{2}\right)
\left(  r^{2}-z_{d}^{2}\right) \nonumber\\
&  =y^{2}+z_{d}^{2}+\beta^{2}r^{2}+2yr+2\beta yz_{d}+2\beta rz_{d}.
\label{Imrep1}%
\end{align}
For the real part, we obtain%
\[
\operatorname{Re}\chi=y^{2}+\operatorname{Re}\left(  z_{d}^{2}\right)
+\beta^{2}\operatorname{Re}\left(  r^{2}\right)  +2y\operatorname{Re}r+2\beta
y\operatorname{Re}z_{d}+2\beta\operatorname{Re}\left(  rz_{d}\right)  .
\]

Let $\left(  \mathbf{x},\mathbf{y}\right)  \in\overrightarrow{\mathcal{E}%
}_{\mathbf{a},\mathbf{b}}\times\left[  \mathbf{c},\mathbf{d}\right]  $,
$\mathbf{z}:=\mathbf{x-Ry}$, and $r_{\operatorname*{R}}:=\left\Vert
\mathbf{z}\right\Vert $. As a first condition we impose $C_{\mathcal{E}}>3$;
then a straightforward calculation shows that $\frac{\beta^{2}}{C_{\mathcal{E}%
}\left(  1+\beta\right)  ^{3}}<\min\left\{  1/6,1/\left(  3\beta\right)
\right\}  $ and all estimates in Lemmata \ref{LemExtNorm} and \ref{Lemrrplus}
are applicable. This implies the existence\footnote{We omit a sharper
specification of the involved constants for the sake of readability.} of some
$C_{1}\geq1$ such that%
\begin{equation}%
\begin{array}
[c]{ll}%
\left\vert \operatorname{Im}r\right\vert  & \overset{\text{(\ref{imr1})}}%
{\leq}C_{1}\kappa r_{\operatorname*{R}},\\
\left\vert \operatorname{Im}z_{d}\right\vert  & \overset{\text{(\ref{estimzd}%
)}}{\leq}C_{1}\kappa r_{\operatorname*{R}},\\
\operatorname{Re}r & \overset{\text{(\ref{rests1})}}{\leq}C_{1}%
r_{\operatorname*{R}},\\
\operatorname{Re}r & \overset{\text{(\ref{rests1})}}{\geq}C_{1}^{-1}%
r_{\operatorname*{R}},\\
\operatorname{Re}\left(  r^{2}\right)  & \overset{\text{(\ref{ReImEsta})}%
}{\geq}\left(  C_{1}^{-1}r_{\operatorname*{R}}\right)  ^{2},\\
\operatorname{Re}z_{d} & \overset{\text{(\ref{realpartfromabove})}}{\leq}%
C_{1}r_{\operatorname*{R}},\\
\operatorname{Re}z_{d} & =\xi_{d}+y_{d}+\operatorname{Re}x_{d}-\xi_{d}%
\overset{\xi_{d}+y_{d}>0}{>}-\left\vert x_{d}-\xi_{d}\right\vert
\overset{\text{(\ref{xidfinal})}}{\geq}-C_{1}\kappa r_{\operatorname*{R}},\\
\operatorname{Re}\left(  rz_{d}\right)  & =\operatorname{Re}r\operatorname{Re}%
z_{d}-\operatorname{Im}r\operatorname{Im}z_{d}\geq-C_{1}^{2}\kappa
r_{\operatorname*{R}}^{2},\\
\left\vert \operatorname{Im}\left(  rz_{d}\right)  \right\vert  &
\leq\left\vert \operatorname{Im}r\right\vert \left\vert \operatorname{Re}%
z_{d}\right\vert +\left\vert \operatorname{Re}r\right\vert \left\vert
\operatorname{Im}z_{d}\right\vert \leq C_{1}^{2}\kappa r_{\operatorname*{R}%
}^{2}.
\end{array}
\label{fullests}%
\end{equation}
In this way, we obtain%
\begin{align}
\operatorname{Re}\chi &  \geq y^{2}-\left(  C_{1}\kappa r_{\operatorname*{R}%
}\right)  ^{2}+\beta^{2}\left(  C_{1}^{-1}r_{\operatorname*{R}}\right)
^{2}+2C_{1}^{-1}yr_{\operatorname*{R}}+2\beta y\left(  -C_{1}\kappa
r_{\operatorname*{R}}\right)  +2\beta\left(  -C_{1}^{2}\kappa
r_{\operatorname*{R}}^{2}\right) \nonumber\\
&  \geq y^{2}+\left(  \left(  \frac{\beta}{C_{1}}\right)  ^{2}-C_{1}^{2}%
\kappa\left(  2\beta+\kappa\right)  \right)  r_{\operatorname*{R}}%
^{2}+2yr_{\operatorname*{R}}\left(  C_{1}^{-1}-\beta C_{1}\kappa\right)  .
\label{Rechilow}%
\end{align}
For the prefactors it holds for the considered range of $\kappa:$%
\begin{align}
\left(  \frac{\beta}{C_{1}}\right)  ^{2}-C_{1}^{2}\kappa\left(  2\beta
+\kappa\right)   &  \geq\beta^{2}\left(  C_{1}^{-2}-C_{1}^{2}\frac{\beta
}{C_{\mathcal{E}}\left(  1+\beta\right)  ^{3}}\left(  2+\frac{\beta
}{C_{\mathcal{E}}\left(  1+\beta\right)  ^{3}}\right)  \right) \nonumber\\
&  \geq\beta^{2}\left(  C_{1}^{-2}-\frac{C_{1}^{2}}{C_{\mathcal{E}}}\left(
2+\frac{1}{C_{\mathcal{E}}}\right)  \right) \label{betaC1}\\
C_{1}^{-1}-\beta C_{1}\kappa &  \geq C_{1}^{-1}-\frac{C_{1}}{C_{\mathcal{E}}}.
\label{betaC2}%
\end{align}
A possible adjustment of $C_{\mathcal{E}}$ (only depending on the number
$C_{1}>0$) implies that the right-hand side in (\ref{betaC2}) is positive (so
that the last term in (\ref{Rechilow}) can be dropped in a lower bound) and
there is a number $0<C_{2}\leq1$ such that the right-hand side in
(\ref{betaC1}) is bounded from below by $C_{2}\beta^{2}$. In this way%
\begin{equation}
\operatorname{Re}\chi\geq y^{2}+C_{2}\beta^{2}r_{\operatorname*{R}}^{2}
\label{R2chi2}%
\end{equation}
follows.

For the imaginary part, we start from (\ref{Imrep1}) and obtain%
\[
\operatorname{Im}\chi=2\operatorname{Re}z_{d}\operatorname{Im}z_{d}+2\beta
^{2}\operatorname{Re}r\operatorname{Im}r+2y\operatorname{Im}r+2\beta
y\operatorname{Im}z_{d}+2\beta\operatorname{Im}\left(  rz_{d}\right)  .
\]
This implies the estimate of the modulus%
\begin{align*}
\left\vert \operatorname{Im}\chi\right\vert  &  \leq2C_{1}^{2}\kappa
r_{\operatorname*{R}}^{2}+2\beta^{2}C_{1}^{2}\kappa r_{\operatorname*{R}}%
^{2}+2C_{1}\kappa yr_{\operatorname*{R}}+2C_{1}\beta\kappa
yr_{\operatorname*{R}}+2C_{1}^{2}\beta\kappa r_{\operatorname*{R}}^{2}\\
&  \leq\kappa\left(  2yC_{1}\left(  1+\beta\right)  r_{\operatorname*{R}%
}+2C_{1}^{2}\left(  1+\beta\right)  ^{2}r_{\operatorname*{R}}^{2}\right) \\
&  \leq\kappa\left(  y^{2}+3\left(  1+\beta\right)  ^{2}C_{1}^{2}%
r_{\operatorname*{R}}^{2}\right)  .
\end{align*}
The ratio of the imaginary and real part of $\chi$ can be bounded by%
\begin{equation}
\frac{\left\vert \operatorname{Im}\chi\right\vert }{\operatorname{Re}\chi}%
\leq\frac{y^{2}+3\left(  1+\beta\right)  ^{2}C_{1}^{2}r_{\operatorname*{R}%
}^{2}}{y^{2}+C_{2}\beta^{2}r_{\operatorname*{R}}^{2}}\kappa\leq3\frac
{C_{1}^{2}}{C_{2}}\frac{\left(  1+\beta\right)  ^{2}}{\beta^{2}}\kappa.
\label{IMdivRe}%
\end{equation}
For the modulus it holds%
\begin{align*}
\left\vert \chi\right\vert  &  \leq\left\vert y+r_{+}\right\vert ^{2}+\left(
1+\beta^{2}\right)  r_{\operatorname*{R}}^{2}\leq2y^{2}+2\left\vert
r_{+}\right\vert ^{2}+\left(  1+\beta^{2}\right)  r_{\operatorname*{R}}^{2}\\
&  \overset{\text{(\ref{bothboundsnplus}), (\ref{rests1})}}{\leq}%
2y^{2}+3\left(  1+\beta\right)  ^{2}r_{\operatorname*{R}}^{2}.
\end{align*}
We have collected all ingredients for deriving estimates for $\sqrt{\chi}$.
For the modulus we obtain%
\begin{equation}
\left\vert \sqrt{\chi}\right\vert \overset{\text{(\ref{normz2x})}}{=}%
\sqrt{\left\vert \chi\right\vert }\leq\sqrt{3}\left(  y+\left(  1+\beta
\right)  r_{\operatorname*{R}}\right)  \label{estmodSqrtchi}%
\end{equation}
and for the real part%
\begin{align}
\operatorname{Re}\sqrt{\chi}  &  \overset{\text{(\ref{Reestfrombelow})}}{\geq
}\sqrt{\left\vert \chi\right\vert }\left(  1-\frac{1}{8}\left(  \frac
{\operatorname{Im}\chi}{\operatorname{Re}\chi}\right)  ^{2}\right)  \geq
\sqrt{\operatorname{Re}\chi}\left(  1-\frac{1}{8}\left(  \frac
{\operatorname{Im}\chi}{\operatorname{Re}\chi}\right)  ^{2}\right) \nonumber\\
&  \overset{\text{(\ref{R2chi2}), (\ref{IMdivRe})}}{\geq}\left(  1-\frac{1}%
{8}\left(  3\frac{C_{1}^{2}}{C_{2}}\frac{\left(  1+\beta\right)  ^{2}}%
{\beta^{2}}\kappa\right)  ^{2}\right)  \sqrt{y^{2}+C_{2}\beta^{2}%
r_{\operatorname*{R}}^{2}}. \label{estReSqrtmue}%
\end{align}
For the considered range of $\kappa$ it holds%
\[
\frac{1}{8}\left(  3\frac{C_{1}^{2}}{C_{2}}\frac{\left(  1+\beta\right)  ^{2}%
}{\beta^{2}}\kappa\right)  ^{2}\leq\frac{9}{8}\frac{C_{1}^{4}}{C_{\mathcal{E}%
}^{2}C_{2}^{2}}.
\]
Again, by a possible adjustment of $C_{\mathcal{E}}$ only depending on the
numbers $C_{1}$, $C_{2}$, the prefactor in (\ref{estReSqrtmue}) is bounded
from below by $1/2$ and%
\begin{equation}
\operatorname{Re}\sqrt{\chi}\geq\frac{1}{2}\sqrt{y^{2}+C_{2}\beta
^{2}r_{\operatorname*{R}}^{2}}. \label{ReSqrtchifrombelow}%
\end{equation}
Using the definitions in (\ref{muecomprep}) it is a straightforward
calculation to verify $t+\beta\tilde{\mu}=\sqrt{\chi}$ and (\ref{mueestb2}) follows.

Next the modulus $\left\vert \tilde{\mu}\left(  \mathbf{z},t\right)
\right\vert $ will be estimated and the representation of $\tilde{\mu}$ in
(\ref{muecomprep}) is employed. For the numerator, we get%
\begin{align*}
\left\vert a\sqrt{\chi}+\beta\omega^{2}\right\vert  &  \overset
{\text{(\ref{estmodSqrtchi})}}{\leq}\sqrt{3}\left(  y+\left\vert
r_{+}\right\vert \right)  \left(  y+\left(  1+\beta\right)
r_{\operatorname*{R}}\right)  +\beta r_{\operatorname*{R}}^{2}\\
&  \overset{\text{(\ref{bothboundsnplus}), (\ref{rests1})}}{\leq}\sqrt
{3}\left(  y+\left(  1+\beta\right)  r_{\operatorname*{R}}\right)  ^{2}+\beta
r_{\operatorname*{R}}^{2}\\
&  \leq3\left(  y+\left(  1+\beta\right)  r_{\operatorname*{R}}\right)  ^{2}%
\end{align*}
and for the denominator:%
\begin{align*}
\left\vert \beta a+\sqrt{\chi}\right\vert  &  \geq\beta\left(
y+\operatorname{Re}r_{+}\right)  +\operatorname{Re}\sqrt{\chi}\\
&  \overset{\text{(\ref{fullests}), (\ref{ReSqrtchifrombelow})}}{\geq}%
\beta\left(  y+\frac{1-\beta C_{1}^{2}\kappa}{C_{1}}r_{\operatorname*{R}%
}\right)  +\frac{1}{2}\sqrt{y^{2}+C_{2}\beta^{2}r_{\operatorname*{R}}^{2}}.
\end{align*}
The bound on $\kappa$ implies $C_{1}\beta\kappa\leq C_{1}/C_{\mathcal{E}}$ and
for sufficiently large $C_{\mathcal{E}}$ (only depending on $C_{1}$) it holds
$1-\beta C_{1}^{2}\kappa\geq1/2$ so that%
\begin{equation}
\operatorname{Re}r_{+}\overset{\text{(\ref{fullests})}}{\geq}\left(
C_{1}^{-1}-C_{1}\beta\kappa\right)  r_{\operatorname*{R}}\geq\frac
{r_{\operatorname*{R}}}{2C_{1}}. \label{estrplusbelowser}%
\end{equation}
Thus,%
\begin{align}
\left\vert \beta a+\sqrt{\chi}\right\vert  &  \geq\beta\left(  y+\frac
{r_{\operatorname*{R}}}{2C_{1}}\right)  +\frac{1}{2\sqrt{2}}\left(
y+\sqrt{C_{2}}\beta r_{\operatorname*{R}}\right) \label{rpluslow}\\
&  \geq\left(  \frac{1}{2\sqrt{2}}+\beta\right)  y+\beta\left(  \frac
{\sqrt{C_{2}}}{2\sqrt{2}}+\frac{1}{2C_{1}}\right)  r_{\operatorname*{R}%
}\nonumber\\
&  \geq C_{3}\left(  \left(  1+\beta\right)  y+\beta r_{\operatorname*{R}%
}\right) \nonumber
\end{align}
for some $C_{3}$ depending only on $C_{1}$ and $C_{2}$. This leads to the
upper estimate
\[
\left\vert \tilde{\mu}\left(  \mathbf{z},t\right)  \right\vert \leq\frac
{3}{C_{3}\left(  1+\beta\right)  }\frac{y+\left(  1+\beta\right)
r_{\operatorname*{R}}}{y+\frac{\beta}{\left(  1+\beta\right)  }%
r_{\operatorname*{R}}}\left(  y+\left(  1+\beta\right)  r_{\operatorname*{R}%
}\right)  \leq\frac{3}{C_{3}}\frac{1+\beta}{\beta}\left(  y+\left(
1+\beta\right)  r_{\operatorname*{R}}\right)  .
\]

To derive a lower estimate for the real part we use%
\begin{equation}
\operatorname{Re}\tilde{\mu}\left(  \mathbf{z},t\right)  =\frac
{\operatorname{Re}\left(  \left(  a\sqrt{\chi}+\beta\omega^{2}\right)  \left(
\overline{\beta a+\sqrt{\chi}}\right)  \right)  }{\left\vert \beta
a+\sqrt{\chi}\right\vert ^{2}}=\frac{N}{\left\vert \beta a+\sqrt{\chi
}\right\vert ^{2}} \label{EqRemuetilde}%
\end{equation}
with%
\[
N:=\beta\left\vert a\right\vert ^{2}\operatorname{Re}\sqrt{\chi}+\left\vert
\chi\right\vert \operatorname{Re}a+\beta^{2}\operatorname{Re}\left(
a\omega^{2}\right)  +\beta\operatorname{Re}\left(  \omega^{2}\overline
{\sqrt{\chi}}\right)
\]
and estimate the terms in $N$ step by step. From (\ref{fullests}) one
concludes that%
\begin{align}
\left\vert a\right\vert ^{2}  &  \geq\left(  y+\operatorname{Re}r_{+}\right)
^{2}\overset{\text{(\ref{estrplusbelowser})}}{\geq}\left(  y+\frac
{r_{\operatorname*{R}}}{2C_{1}}\right)  ^{2},\nonumber\\
\operatorname{Re}\sqrt{\chi}  &  \overset{\text{(\ref{ReSqrtchifrombelow})}%
}{\geq}\frac{1}{2}\sqrt{y^{2}+C_{2}\beta^{2}r_{\operatorname*{R}}^{2}}%
\geq\frac{1}{3}\left(  y+\sqrt{C_{2}}\beta r_{\operatorname*{R}}\right)
,\label{EstReSqrtchi}\\
\left\vert \chi\right\vert  &  \geq\operatorname{Re}\chi\overset
{\text{(\ref{R2chi2})}}{\geq}y^{2}+C_{2}\beta^{2}r_{\operatorname*{R}}%
^{2},\nonumber\\
\operatorname{Re}a  &  \geq y+\frac{r_{\operatorname*{R}}}{2C_{1}}.
\label{EstRea}%
\end{align}
This allows us to estimate the first two summands in the definition of $N$%
\begin{align*}
\beta\left\vert a\right\vert ^{2}\operatorname{Re}\sqrt{\chi}+\left\vert
\chi\right\vert \operatorname{Re}a  &  \geq\left(  y+\frac
{r_{\operatorname*{R}}}{2C_{1}}\right)  \left(  \frac{\beta}{3}\left(
y+\frac{r_{\operatorname*{R}}}{2C_{1}}\right)  \left(  y+\sqrt{C_{2}}\beta
r_{\operatorname*{R}}\right)  +\left(  y^{2}+C_{2}\beta^{2}%
r_{\operatorname*{R}}^{2}\right)  \right) \\
&  \geq\left(  y+\frac{r_{\operatorname*{R}}}{2C_{1}}\right)  \left(
\frac{\beta}{3}\left(  y^{2}+\frac{\sqrt{C_{2}}}{2C_{1}}\beta
r_{\operatorname*{R}}^{2}\right)  +y^{2}+C_{2}\beta^{2}r_{\operatorname*{R}%
}^{2}\right) \\
&  \geq\left(  y+\frac{r_{\operatorname*{R}}}{2C_{1}}\right)  \left(  \left(
1+\frac{\beta}{3}\right)  y^{2}+\left(  \frac{\sqrt{C_{2}}}{6C_{1}}%
+C_{2}\right)  \beta^{2}r_{\operatorname*{R}}^{2}\right) \\
&  \geq C_{4}\left(  y+\frac{r_{\operatorname*{R}}}{2C_{1}}\right)  \left(
\left(  1+\beta\right)  y^{2}+\beta^{2}r_{\operatorname*{R}}^{2}\right)
\end{align*}
for some $C_{4}$ depending only on $C_{1}$ and $C_{2}$. To estimate the last
two summands in $N$ we start with the relations%
\begin{align*}
\operatorname{Re}\left(  a\omega^{2}\right)   &  =\operatorname{Re}%
a\operatorname{Re}\left(  \omega^{2}\right)  -\operatorname{Im}%
a\operatorname{Im}\left(  \omega^{2}\right) \\
\operatorname{Re}\left(  \omega^{2}\overline{\sqrt{\chi}}\right)   &
=\operatorname{Re}\left(  \omega^{2}\right)  \operatorname{Re}\sqrt{\chi
}+\operatorname{Im}\sqrt{\chi}\operatorname{Im}\left(  \omega^{2}\right)  .
\end{align*}
For the single factors, we employ (\ref{ReImEst}) and%
\begin{align*}
\operatorname{Re}\left(  z_{d}^{2}\right)   &  =\left(  \operatorname{Re}%
z_{d}\right)  ^{2}-\left(  \operatorname{Im}z_{d}\right)  ^{2}\overset
{\text{(\ref{realpartfromabove})}}{\leq}r_{\operatorname*{R}}^{2},\\
\left\vert \operatorname{Re}a\right\vert  &  \leq y+\operatorname{Re}\left(
r+\beta z_{d}\right)  \overset{\text{(\ref{rests1}), (\ref{realpartfromabove}%
)}}{\leq}y+\left(  1+\beta\right)  r_{\operatorname*{R}}\\
\left\vert \operatorname{Im}a\right\vert  &  \leq\operatorname{Im}\left(
r+\beta z_{d}\right)  \overset{\text{(\ref{imr1}), (\ref{estimzd})}}{\leq
}2\kappa\left(  2+\beta\right)  r_{\operatorname*{R}}%
\end{align*}
to get%
\[
\operatorname{Re}\left(  \omega^{2}\right)  \geq\operatorname{Re}\left\langle
\mathbf{z},\mathbf{z}\right\rangle -\operatorname{Re}\left(  z_{d}^{2}\right)
\geq\left(  1-8\kappa^{2}\right)  r_{\operatorname*{R}}^{2}%
-r_{\operatorname*{R}}^{2}\geq-8\kappa^{2}r_{\operatorname*{R}}^{2}%
\quad\text{and\quad}\left\vert \operatorname{Im}\left(  \omega^{2}\right)
\right\vert \leq4\kappa r_{\operatorname*{R}}^{2},
\]
while the estimate of $\sqrt{\chi}$ follows in a similar fashion as
(\ref{estReSqrtmue}) from (\ref{estimpart})%
\[
\left\vert \operatorname{Im}\sqrt{\chi}\right\vert \leq\frac{\left\vert
\operatorname{Im}\chi\right\vert }{2\operatorname{Re}\chi}\sqrt{\left\vert
\chi\right\vert }\overset{\text{(\ref{IMdivRe}), (\ref{estmodSqrtchi})}}{\leq
}\frac{3\sqrt{3}C_{1}^{2}}{2C_{2}}\frac{\left(  1+\beta\right)  ^{2}}%
{\beta^{2}}\kappa\left(  y+\left(  1+\beta\right)  r_{\operatorname*{R}%
}\right)  .
\]
In this way, we obtain for the real parts of the products:%
\begin{align*}
\operatorname{Re}\left(  a\omega^{2}\right)   &  \overset
{\text{(\ref{bothboundsnplus}), (\ref{fullests})}}{\geq}-C_{5}\left(
y+\left(  1+\beta\right)  r_{\operatorname*{R}}\right)  \kappa^{2}%
r_{\operatorname*{R}}^{2},\\
\operatorname{Re}\left(  \omega^{2}\overline{\sqrt{\chi}}\right)   &
\overset{\text{(\ref{estmodSqrtchi})}}{\geq}-C_{5}\frac{\left(  1+\beta
\right)  ^{2}}{\beta^{2}}\left(  y+\left(  1+\beta\right)
r_{\operatorname*{R}}\right)  \kappa^{2}r_{\operatorname*{R}}^{2}%
\end{align*}
for a constant $C_{5}$ which only depends on the numbers $C_{1}$ and $C_{2}$.
The sum of both terms can be estimated by%
\begin{align*}
\operatorname{Re}\left(  a\omega^{2}\right)  +\operatorname{Re}\left(
\omega^{2}\overline{\sqrt{\chi}}\right)   &  \geq-2C_{5}\frac{\left(
1+\beta\right)  ^{2}}{\beta^{2}}\left(  y+\left(  1+\beta\right)
r_{\operatorname*{R}}\right)  \kappa^{2}r_{\operatorname*{R}}^{2}\\
&  \geq-\frac{2C_{5}}{C_{\mathcal{E}}^{2}}\frac{\beta^{2}}{\left(
1+\beta\right)  ^{4}}\left(  y+\left(  1+\beta\right)  r_{\operatorname*{R}%
}\right)  r_{\operatorname*{R}}^{2}.
\end{align*}
Thus, for the numerator $N$ in (\ref{EqRemuetilde}) it follows%
\begin{align*}
N  &  \geq C_{4}\left(  y+\frac{r_{\operatorname*{R}}}{2C_{1}}\right)  \left(
\left(  1+\beta\right)  y^{2}+\beta^{2}r_{\operatorname*{R}}^{2}\right)
-\frac{2C_{5}}{C_{\mathcal{E}}^{2}}\frac{1}{\left(  1+\beta\right)  ^{3}%
}\left(  y+r_{\operatorname*{R}}\right)  \beta^{2}r_{\operatorname*{R}}^{2}\\
&  \geq\left(  C_{4}\left(  y+\frac{r_{\operatorname*{R}}}{2C_{1}}\right)
-\frac{2C_{5}}{C_{\mathcal{E}}}\left(  y+\frac{r_{\operatorname*{R}}%
}{C_{\mathcal{E}}}\right)  \right)  \left(  \left(  1+\beta\right)
y^{2}+\beta^{2}r_{\operatorname*{R}}^{2}\right)  .
\end{align*}
Again, by a possible adjustment of $C_{\mathcal{E}}$ only depending on the
numbers $C_{1}$, $C_{2}$, $C_{4}$, $C_{5}$ we end up with%
\[
N\geq C_{6}\left(  y+r_{\operatorname*{R}}\right)  \left(  \left(
1+\beta\right)  y^{2}+\beta^{2}r_{\operatorname*{R}}^{2}\right)  .
\]
For the denominator it holds%
\begin{equation}
\left\vert \beta a+\sqrt{\chi}\right\vert \leq\beta\left(  y+\left\vert
r_{+}\right\vert \right)  +\left\vert \sqrt{\chi}\right\vert \overset
{\text{(\ref{rplusupperbound}), (\ref{estmodSqrtchi})}}{\leq}\left(
\beta+\sqrt{3}\right)  \left(  y+\left(  1+\beta\right)  r_{\operatorname*{R}%
}\right)  . \label{betaplussqrtchiupper}%
\end{equation}
The combination of these two inequalities leads to%
\begin{align*}
\operatorname{Re}\tilde{\mu}\left(  \mathbf{z},t\right)   &  \geq\frac
{C_{6}\left(  y+r_{\operatorname*{R}}\right)  \left(  \left(  1+\beta\right)
y^{2}+\beta^{2}r_{\operatorname*{R}}^{2}\right)  }{\left(  \beta+\sqrt
{3}\right)  ^{2}\left(  y+\left(  1+\beta\right)  r_{\operatorname*{R}%
}\right)  ^{2}}\geq\frac{\left(  1+\beta\right)  C_{6}}{2\left(  \beta
+\sqrt{3}\right)  ^{2}}\left(  \frac{y+\frac{\beta}{\sqrt{1+\beta}%
}r_{\operatorname*{R}}}{y+\left(  1+\beta\right)  r_{\operatorname*{R}}%
}\right)  ^{2}\left(  y+r_{\operatorname*{R}}\right) \\
&  \geq\frac{C_{6}}{2\left(  \beta+\sqrt{3}\right)  ^{2}}\frac{\beta^{2}%
}{\left(  1+\beta\right)  ^{2}}\left(  y+r_{\operatorname*{R}}\right)
\end{align*}
and the assertion for $\operatorname{Re}\tilde{\mu}$ follows.\bigskip

Next, we estimate the modulus of $\tilde{\mu}^{\prime}$ and recall
(\ref{defsmue}):%
\begin{align}
\tilde{\mu}^{\prime}\left(  \mathbf{z},y\right)   &  :=\frac{\partial
\tilde{\mu}\left(  \mathbf{z},y\right)  }{\partial y}=\frac{t\left(
\mathbf{z},y\right)  }{t\left(  \mathbf{z},y\right)  +\beta\tilde{\mu}\left(
\mathbf{z},y\right)  }\nonumber\\
&  \overset{\text{(\ref{defsmue}), (\ref{muecomprep})}}{=}\frac{a^{2}%
-\omega^{2}}{a\left(  a+\beta\sqrt{\chi}\right)  +\left(  \beta^{2}-1\right)
\omega^{2}}=\frac{a^{2}-\omega^{2}}{\sqrt{\chi}\left(  \beta a+\sqrt{\chi
}\right)  }. \label{mueprimerep2}%
\end{align}
For the numerator, we get%
\[
\left\vert a^{2}-\omega^{2}\right\vert \leq\left\vert a\right\vert
^{2}+\left\vert \omega\right\vert ^{2}\overset{\text{(\ref{rplusupperbound})}%
}{\leq}\left(  y+\left(  1+\beta\right)  r_{\operatorname*{R}}\right)
^{2}+r_{\operatorname*{R}}^{2}\leq2\left(  y+\left(  1+\beta\right)
r_{\operatorname*{R}}\right)  ^{2}%
\]
and for the denominator:%
\[
\left\vert \sqrt{\chi}\right\vert \left\vert \beta a+\sqrt{\chi}\right\vert
\overset{\text{(\ref{EstReSqrtchi}),(\ref{rpluslow})}}{\geq}\frac{C_{3}}%
{3}\left(  y+\sqrt{C_{2}}\beta r_{\operatorname*{R}}\right)  \left(  \left(
1+\beta\right)  y+\beta r_{\operatorname*{R}}\right)  .
\]
This leads to%
\begin{align*}
\left\vert \tilde{\mu}^{\prime}\left(  \mathbf{z},y\right)  \right\vert  &
\leq\frac{6}{C_{3}}\frac{\left(  y+\left(  1+\beta\right)
r_{\operatorname*{R}}\right)  ^{2}}{\left(  y+\sqrt{C_{2}}\beta
r_{\operatorname*{R}}\right)  \left(  \left(  1+\beta\right)  y+\beta
r_{\operatorname*{R}}\right)  }\leq\frac{6}{C_{3}}\frac{1+\beta}{\beta}%
\frac{y+\left(  1+\beta\right)  r_{\operatorname*{R}}}{y+\sqrt{C_{2}}\beta
r_{\operatorname*{R}}}\\
&  \leq\frac{6}{C_{3}}\frac{1+\beta}{\beta}\max\left\{  1,\frac{1+\beta}%
{\sqrt{C_{2}}\beta}\right\}
\end{align*}
and finally to the assertion
\endproof

\bibliographystyle{abbrv}
\bibliography{nlailu}

\end{document}